\documentclass[11pt]{article}
\usepackage{amsmath}
\usepackage{amsmath,amssymb,latexsym,color}
\usepackage[mathscr]{eucal}
\newtheorem{thm}{Theorem}[section]

\newtheorem{lemma}[thm]{Lemma}

\newtheorem{example}{Example}[section]
\newtheorem{defin}{Definition}[section]

\newtheorem{remark}{Remark}[section]

\newcommand{\proof}{{\it Proof.\quad}}
\newcommand{\qed}{\hfill\Box\medskip}

\usepackage{CJK}
\begin{document}
\begin{CJK*}{GBK}{song}

\newcommand{\be}{\begin{equation}\label}
\newcommand{\ee}{\end{equation}}
\newcommand{\bea}{\begin{eqnarray}\label}
\newcommand{\eea}{\end{eqnarray}}

\title{\bf The smallest one-realization of a given set \uppercase\expandafter{\romannumeral4} }

\author{
Kefeng Diao$^{\rm a}$\quad Vitaly I. Voloshin$^{\rm b}$\quad Kaishun
Wang$^{\rm c}$\thanks{Corresponding
author.\newline    Email address: kfdiao@163.com, vvoloshin@troy.edu, wangks@bnu.edu.cn, zhaopingly@163.com}\quad Ping Zhao$^{\rm a}$\\
{\footnotesize a. \em School of Science, Linyi  University,
Linyi, Shandong, 276005, China }\\
{\footnotesize b. \em  Department of Mathematics and Geomatics, Troy
University, Troy, AL 36082, U.S.A.}\\
{\footnotesize c. \em  Sch. Math. Sci. {\rm \&} Lab. Math. Com.
Sys., Beijing Normal University, Beijing 100875,  China} }

\date{}
 \maketitle

\begin{abstract}

Let  $S$ be a finite set of positive integers. A mixed hypergraph
${\cal H}$ is a one-realization of $S$ if its feasible set is $S$
and each entry of its chromatic spectrum is either 0 or 1. In [P.
Zhao, K. Diao, Y. Chang and K. Wang, The smallest one-realization of
a given set \uppercase\expandafter{\romannumeral2},  Discrete Math.
312 (2012) 2946--2951], we determined the minimum number
 of vertices of a $3$-uniform bi-hypergraph  which is a one-realization of $S$.
 In this paper, we generalize this result to $r$-uniform
 bi-hypergraphs.

\medskip
\noindent {\em Key words:} hypergraph coloring; mixed hypergraph; $r$-uniform
bi-hypergraph; feasible set; chromatic spectrum; one-realization
\end{abstract}

\section{Introduction}

 A {\em mixed hypergraph } on a finite set $X$ is a triple ${\cal
H}=(X, {\cal C}, {\cal D})$, where ${\cal C}$ and ${\cal D}$ are
families of subsets of $X$. The members of ${\cal C}$ and ${\cal D}$
are called {\em ${\cal C}$-edges} and {\em ${\cal D}$-edges},
respectively, and all have the cardinality at least 2. For   a subset $X'$ of $X$, a hypergraph ${\cal
H}'=(X', {\cal C}', {\cal D}')$  is called an {\em induced
sub-hypergraph}
  of  ${\cal H}$ on $X'$, denoted by ${\cal H}[X']$, if ${\cal C}'$ (resp. ${\cal
  D}'$) consists of   ${\cal C}$-edges (resp. ${\cal D}$-edges)
  of $\cal H$ whose vertices belong to $X'.$    A mixed hypergraph
 is \emph{$r$-uniform} if every
edge contains $r$ vertices. A mixed hypergraph ${\cal H}=(X, {\cal
B}, {\cal B})$ is called a \emph{bi-hypergraph}, denoted by ${\cal
H}=(X, {\cal B})$, whose edges are called \emph{bi-edges}.

   The distinction between ${\cal C}$-edges and ${\cal
D}$-edges becomes substantial when colorings
 are considered. A {\em proper $k$-coloring} of
${\cal H}$ is a mapping from $X$ into a set of $k$ colors such that
each ${\cal C}$-edge has at least two vertices with a \emph{Common}
color and each ${\cal D}$-edge has at least two vertices with
\emph{Distinct} colors.  A coloring may also be viewed as a
\emph{partition} (\emph{feasible partition}) of $X$, where the
\emph{color classes} (partition classes) are the sets of vertices
assigned to the same color. A {\em strict $k$-coloring} is a proper
$k$-coloring with $k$ nonempty color classes, and a mixed hypergraph
is {\em $k$-colorable} if it has a strict $k$-coloring. An edge  is said to be \emph{monochromatic} (resp. \emph{polychromatic}) if all its vertices have the same color (resp. different colors).
  The set of all the  values $k$ such that ${\cal H}$ has a strict $k$-coloring is called the {\em feasible set}
  of ${\cal H}$, denoted by ${\cal F}({\cal H})$. The {\em lower chromatic number} $\chi ({\cal H})$ of ${\cal H}$ is the minimum number contained in ${\cal F}({\cal H})$ and the  {\em upper chromatic number} $\overline \chi ({\cal H})$ of ${\cal H}$ is the maximum number contained in ${\cal F}({\cal H})$. For each $k$, let $r_k$ denote
the number of {\em feasible partitions} of the vertex set of  ${\cal H}$. The
vector $R({\cal H})=(r_1,r_2,\ldots,r_{\overline\chi})$ is called
the {\em chromatic spectrum} of ${\cal H}$, where $\overline\chi$ is the  upper chromatic number of ${\cal H}$. If $S$ is a finite set of positive integers, we say that a
mixed hypergraph ${\cal H}$ is a {\it realization} of $S$ if ${\cal
F}({\cal H})=S$. It is readily seen that if $1 \in {\cal
F}(\cal{H})$, then $\cal{H}$ cannot have any $\cal{D}$-edges. A
mixed hypergraph ${\cal H}$ is a {\it one-realization} of $S$ if it
is a realization of $S$ and each entry of the chromatic
spectrum of ${\cal H}$ is either 0 or 1. The study of   the
colorings of mixed hypergraphs has made a lot of progress since its
inception \cite{Voloshin1}. The applications of mixed hypergraph coloring include modeling of
several types of graph coloring (like list coloring without lists), different kinds of homomorphism of simple graphs and multigraphs,
channel assignment problems \cite{Kral2007} and the newest application in problems arising in cyber security \cite{Jaffe}.
For more information, see \cite{Voloshin2} and the regularly updated website \cite{Volweb}.

K\"{u}ndgen et al. \cite{Kundgen} initiated  a systematic study
of planar mixed hypergraphs, and found a one-realization of
$\{2,4\}$ on 6 vertices for planar hypergraphs. Let $S$ be a finite
set of positive integers with $\min(S)\geq 2.$
 Bacs$\acute{\rm o}$
et al. \cite{Bacso} discussed the properties of uniform
bi-hypergraphs
       ${\cal H}$ which are one-realizations of $S$ when $|S|=1$, in this case we also say that
       ${\cal H}$ is \emph{uniquely colorable}.
       Jiang et al. \cite{Jiang} proved that the minimum number of
  vertices of realizations  of
$S$  is $2\max(S)-\min(S)$ if $|S|=2$ and $\max(S)-1\notin S$.
Kr\'{a}l \cite{Kral} proved that there exists a
one-realization of $S$ with at most $|S|+2\max(S)-\min(S)$ vertices.
In \cite{zdw2} we   proved that the minimum number of vertices of
one-realizations of $S$ is $2\max(S)-\min(S)$ if $\max(S)-1\notin
 S$ or $2\max(S)-\min(S)-1$ if $\max(S)-1\in S$.   Bujt\'{a}s and Tuza
\cite {Bujtas} proved that $S$ is the feasible set of some $r$-uniform bi-hypergraph
${\cal H}=(X, {\cal B})$ with ${\cal B}\neq \emptyset$ if and only if
(i) $\min(S)\geq r$, or (ii) $2\leq \min(S)\leq r-1$ and $S$ contains all integers between $\min(S)$ and $r-1$.
They also raised the following open problem:

\textbf{Problem.} Let $r$ be an integer at least $3$. Determine the minimum number of vertices
in an $r$-uniform bi-hypergraph with a given feasible set.

The motivation of this paper is to solve the above problem.
 In \cite{zdw1}, we constructed a family of 3-uniform bi-hypergraphs with a given chromatic spectrum,
 and obtained
an upper bound on the minimum number of vertices of a
one-realization  of a given set. Recently, in \cite{zdw3}  we
determined the minimum number of vertices of a $3$-uniform
bi-hypergraph  which is a one-realization  of a given set.
 In this paper, we generalize this result to $r$-uniform bi-hypergraphs and get the following result.

\begin{thm}\label{thm}
For integers $s\geq 2, r\geq 4$ and $n_1>n_2>\cdots >n_s$, let
$\delta_r(S)$ be the minimum number of vertices
 of an $r$-uniform bi-hypergraph  which is a one-realization of
 $S:=\{n_1,n_2,\ldots,n_s\}$. Then $n_s\geq r-1$. Moreover,
 $$\delta_r(S)=\left \{
\begin{array}{ll}
(n_1-2)r -(n_1-r-1)\lfloor\frac r2\rfloor+2, & \mbox{if $n_1>n_2+1$,}\\
(n_1-3)r-(n_1-r-3)\lfloor\frac{r}{2}\rfloor+2, & \mbox{if $n_1=n_2+1>r+2$,}\\
r(r-1)-1, & \mbox{if $n_1=r$ or $n_1=n_3+2=r+1>5$,}\\
(n_1-2)(r-1)+3, & \mbox{otherwise.}\\
\end{array}
\right. $$
\end{thm}

This paper is organized as follows. In Section 2, we show the first
statement of Theorem~\ref{thm}, then prove that the number in
Theorem~\ref{thm}  is a lower bound for $\delta_r(S)$.
In Section 3, we introduce a basic construction of $r$-uniform bi-hypergraphs and discuss its coloring property. In Section 4, we construct $r$-uniform bi-hypergraphs which are
one-realizations of $S$ and meet this lower bound in each case.

\section{The lower bound}

 In this section we  show the first statement of Theorem~\ref{thm} and
  the lower bound for $\delta_r(S)$.

\begin{lemma}\label{lem:ineq1} With the notation of Theorem 1.1, we
have $n_s\geq r-1$.
\end{lemma}
\proof   Let ${\cal H}$ be any $r$-uniform bi-hypergraph which is a
one-realization of $S$. Pick a  strict $n_s$-coloring
$c=\{C_1,C_2,\ldots,C_{n_s}\}$  of ${\cal H}$ with $|C_1|\geq
|C_2|\geq\cdots\geq |C_{n_s}|\geq 1$.

 Suppose for a contradiction that $n_s\leq r-2$. Note that $n_s\geq 2$. If $|C_1|\geq 3$, pick  $x_1,x_2,x_3\in C_1$,
 then  $\{\{x_1\}, C_1\setminus \{x_1\},C_2,\ldots,C_{n_s}\}$ and
 $\{\{x_2\}, C_1\setminus \{x_2\},C_2,\ldots,C_{n_s}\}$ are two distinct strict $(n_s+1)$-colorings of ${\cal H}$, a
  contradiction to that ${\cal H}$ is a one-realization of $S$. If $|C_1|=2$, say $C_1=\{x,y\}$,
  then  $\{\{x\}, \{y\}\cup C_2, C_3,\ldots,C_{n_s}\}$ is a strict  $n_s$-coloring of ${\cal H}$ distinct with $c$, also a
  contradiction. If $|C_1|=1$,  then $\{C_1\cup C_2, C_3,\ldots,C_{n_s}\}$ is a strict $(n_s-1)$-coloring of ${\cal H}$,  a
  contradiction to that $n_s-1\notin S$. Hence, the desired result follows.  $\qed$

In the remaining  we   follow the notation of Theorem~\ref{thm}, and
assume that  ${\cal H}:=(X, {\cal B})$ is an $r$-uniform
bi-hypergraph which is a one-realization of  $S$,
$c=\{C_1,\ldots,C_{n_1}\}$ is a strict $n_1$-coloring of ${\cal H}$ with
$|C_1|\geq |C_2|\geq \cdots\geq |C_{n_1}|=t\geq 1$.

\begin{lemma}\label{lem:basic}  Let $C_1'$ be a maximal subset of $X$
 with distinct colors, and let $2\leq i\leq r-1$ be an integer.
 For each  $j=2, \ldots, i$, suppose  $C_j'$ is a maximal  subset of vertices with distinct colors
  in $X\setminus(C_1'\cup\cdots\cup C_{j-1}')$. Then
$$
|X|\geq |C_1'\cup C_2'\cup\cdots\cup C_i'|+(r-1)(r-i-1)+1.
$$
\end{lemma}
\proof Let $C=C_1'\cup C_2'\cup\cdots\cup C_i'$. Suppose for a
contradiction that $|X\setminus C|\leq (r-1)(r-i-1)$. If
$X\setminus C=\emptyset$, then the fact that each $C_j'$ does not contain
any bi-edge of ${\cal H}$ implies that $\{C_1', C_2',\ldots,
C_{i}'\}$ is a strict $i$-coloring of ${\cal H}$. Since $|C_1|\geq
2,$ by exchanging  the two vertices of $C_1\cap (C_1'\cup C_2')$, we
get another strict $i$-coloring of ${\cal H}$, a contradiction to the fact that $\mathcal H$ is a one-realization of $S$. If
$X\setminus C\neq\emptyset$, then $i<r-1$.
Pick a
partition $\{C_{i+1}',\ldots,   C_m'\}$ of $X\setminus C$ such that
$m\leq r-1$ and $1\leq |C_l'|\leq r-1$ for $i+1\leq l\leq m$, then
  $\{C_1', C_2',\ldots,  C_m'\}$ is a strict $m$-coloring of
${\cal H}$. Exchanging  the two vertices of $C_1\cap (C_1'\cup
C_2')$, we get another strict $m$-coloring of ${\cal H}$, also a
contradiction. The proof is completed. $\qed$

 \begin{lemma}\label{lem:ineq}
 The number in Theorem~\ref{thm}  is a lower bound for $\delta_r(S)$ in each case.
\end{lemma}

\proof Write $$\delta_1=(n_1-2)r-(n_1-r-1)\lfloor\frac r2\rfloor+2,~~
 \delta_2=(n_1-3)r-(n_1-r-3)\lfloor\frac{r}{2}\rfloor+2,$$
   $$
   C_j=\{x_{j1},x_{j2},\ldots,x_{js_j}\},~~
   C_i'=\{x_{1i},x_{2i},\ldots,x_{n_1i}\},  j\in [n_1],\; i\in
   [t].$$
  We divide our proof into the following two cases.

 \textbf{Case 1}\,\  $n_1>n_2+1$.

 That is to say, $n_1-1\notin S$. If $s_{n_1-1}<r-t$, then $\{C_1,\ldots,C_{n_1-2},C_{n_1-1}\cup
C_{n_1}\}$ is a strict $(n_1-1)$-coloring of ${\cal
H}$, a contradiction to that $n_1-1\notin S$.  Hence, $s_{n_1-1}\geq r-t$, which implies that
  $s_j\geq r-t$ holds for each  $ j\in [n_1-1].$

If $t\leq \lfloor\frac r2\rfloor$, write $
C_j'=\{x_{1j},x_{2j},\ldots,x_{n_1-1,j}\}$ for  $t+1\leq j\leq r-t,$
then  $C_1'\cup C_2'\cup\cdots\cup  C_{r-t}'$ is of size
$t+(n_1-1)(r-t)$. By Lemma~\ref{lem:basic} one gets  $|X|\geq
t+(n_1-1)(r-t)+(r-1)(t-1)+1\geq \delta_1.$  If
$\lfloor\frac{r}{2}\rfloor<t\leq r-1$, by Lemma~\ref{lem:basic} we
have  $|X| \geq  n_1t+(r-1)(r-t-1)+1 \geq \delta_1$. If $t\geq r$,
then $|X|\geq n_1t>\delta_1.$

\textbf{Case 2}\,\  $n_1=n_2+1$.

If $s_{n_2-1}<r-t$, then $\{C_1,\ldots,C_{n_2-1},C_{n_2}\cup C_{n_1}\}$ and
 $\{C_1,\ldots,C_{n_2-2}, C_{n_2-1}\cup C_{n_1},C_{n_2}\}$ are two strict $n_2$-colorings of ${\cal H}$,
 a contradiction to the fact that $\mathcal H$ is a one-realization of $S$. Hence, $s_{n_2-1}\geq r-t$,
which implies that $s_j\geq r-t$ for any $j\in [n_1-2]$. Let
$s_{n_2}=p$. Then $p\geq t$ and  $s_j\geq r-p$  for each   $j\in [n_1-2].$

 Suppose $t\leq \lfloor\frac r2\rfloor$.
If $p\leq r-t$, write $ C_j'=\{x_{1j},x_{2j},\ldots,x_{n_2,j}\}$ for
$t+1\leq j\leq p$ and $ C_k'=\{x_{1k},x_{2k},\ldots,x_{n_2-1,k}\}$
for $p+1\leq k\leq r-t$. Then  $C_1'\cup C_2'\cup\cdots\cup
C_{r-t}'$ is of size $t+p+(n_1-2)(r-t)$. By Lemma~\ref{lem:basic} we
have  $|X|\geq t+p+(n_1-2)(r-t)+(r-1)(t-1)+1. $
 If $p>r-t$, write   $
C_j'=\{x_{1j},x_{2j},\ldots,x_{n_2,j}\}$ for $t+1\leq j\leq r-t,$ by
Lemma~\ref{lem:basic}  we  have  $|X|\geq  t+(n_1-1)(r-t)+(r-1)(t-1)+1.$ If
$\lfloor\frac{r}{2}\rfloor<t\leq r-1$, by Lemma~\ref{lem:basic} again we
get $ |X|\geq   n_1t+(r-1)(r-t-1)+1.$ If $t\geq r$,
 then $|X|\geq n_1t.$

Moreover, if $n_s=r-1$, we claim that $|X|\geq r(r-1)-1$.
Suppose
 $|X|\leq r(r-1)-2$. Let
$c'=\{S_1,S_2,\ldots,S_{r-1}\}$ be a strict $(r-1)$-coloring of
${\cal H}$ with $1\leq |S_1|\leq |S_2|\leq \cdots \leq |S_{r-1}|$.
Then $|S_1\cup S_2|\leq 2(r-1)$. Pick $x\in S_1, y\in S_2$.
  If  $|S_1|=r-1$, then
    $\{(S_1\setminus \{x\})\cup \{y\},(S_2\setminus \{y\})\cup
    \{x\},S_3,\ldots,S_{r-1}\}$ is a strict $(r-1)$-coloring of
    ${\cal H}$  distinct with $c'$, a contradiction.
 Assume that $|S_1|\leq r-2$. If $|S_2|\geq 2$, then $\{S_1\cup \{y\}, S_2\setminus \{y\},S_3,\ldots,S_{r-1}\}$
  is a strict $(r-1)$-coloring of
    ${\cal H}$  distinct with $c'$, also a contradiction. If $|S_2|=1$, then $\{S_1\cup S_2,S_3,\ldots,S_{r-1}\}$ is a strict $(r-2)$-coloring of ${\cal H}$, a contradiction to that $r-2\notin S$. Therefore, our claim is valid.

     Hence,  if $n_1=n_2+1>r+2$, then $|X|\geq
\delta_2$; if $n_1=r$ or $n_1=n_3+2=r+1>5$, then $|X|\geq
r(r-1)-1$; and otherwise,  $|X|\geq
(n_1-2)(r-1)+3$.

The proof is completed. $\qed$

 \section{The basic construction}

 In this section, we shall construct a family of $r$-uniform bi-hypergraphs and discuss their coloring properties. This construction plays an important role in constructing $r$-uniform bi-hypergraphs which are one-realizations of $S$ and meet the bounds in Lemma~\ref{lem:ineq}.

 We first introduce the construction.
 For any positive integer $n$, let $[n]$ denote the set
$\{1,2,\ldots,n\}$.

 \medskip
 {\bf Construction \uppercase\expandafter{\romannumeral1}.}  Suppose $s\geq 2.$ If $s=2$, we assume that $n_1>n_2+1$. Let $l=\lceil\frac r2\rceil$ and $M_s=[r-2]\cup \{n_s\}$. Write
 \begin{eqnarray*}
V&=&\bigcup\limits_{i\in M_s} V_i^s,~V_i^s=\{(\underbrace{i,\ldots,i}_{s},t)\mid t\in [r-1]\}, i\in M_s,\\
U&=&\bigcup\limits_{j\in [n_s]\setminus M_s} U_j^s,~U_j^s=\{(\underbrace{j,\ldots,j}_s,t)\mid t\in [l]\}, j\in [n_s]\setminus M_s,~\\
W_{pk}^s&=&\{(\underbrace{n_p+k,\ldots,n_p+k}_{p-1},n_p,\ldots,n_s,t)\mid t\in [l]\setminus \{1\}\},\\
&&  p\in [s]\setminus \{1,2\},k\in [n_{p-1}-n_p-1]\cup \{0\},\\
W_{2h}^s&=&\{(n_2+h,n_2,\ldots,n_s,t)\mid t\in [l]\setminus \{1\}\}, h\in [n_1-n_2-2]\cup \{0\},\\
T_1^s&=&\{(n_1-1,n_2,\ldots,n_{s-1},t,t)\mid t=2,\ldots,l-1,n_s\} \hbox{~if~} r>4, \\
T_1^s&=&\{(n_1-1,n_2,\ldots,n_{s-1},n_s,n_s)\} \hbox{~if~} r=4, \\
T_2^s&=&\{(n_1,n_2,\ldots,n_{s-1},t,t)\mid t=l,\ldots,r-2,n_s\},\\
K^s&=&\{(\underbrace{n_p+k,\ldots,n_p+k}_{p-1},\underbrace{1,\ldots,1}_{s+2-p})\mid p\in [s]\setminus \{1\},~k\in [n_{p-1}-n_p-1]\cup \{0\}\},\\
X_{n_1,\ldots,n_s}&=& V\cup U
\cup(\bigcup\limits_{p=3}^s\bigcup\limits_{k=0}^{n_{p-1}-n_p-1}W_{pk}^s)\cup(\bigcup\limits_{h=0}^{n_1-n_2-2}W_{2h}^s)
\cup(\bigcup\limits_{i=1}^2T_i^s)\cup K^s,\\
{\cal B}_{n_1,\ldots,n_s}&=&\{\{\alpha_1,\ldots,\alpha_r\}\mid
 \alpha_q\in X_{n_1,\ldots,n_s}, ~2\leq |\{\alpha_1^j,\ldots,\alpha_r^j\}|\leq r-1, j\in [s+1]\}\\
&\cup& \{\{(n_s,\ldots,n_s,1,1),(1,\ldots,1),
\ldots,(r-2,\ldots,r-2,1),(n_s,\ldots,n_s,1)\}\},
\end{eqnarray*}
 where
$\alpha_q^j$ is the $j$-th entry of the vertex $\alpha_q, q\in [r]$.
 Then ${\cal H}_{n_1,\ldots,n_s}=(X_{n_1,\ldots,n_s}, {\cal
B}_{n_1,\ldots,n_s})$ is an $r$-uniform bi-hypergraph with  $(n_1-2)r -(n_1-r-1)\lfloor\frac r2\rfloor+2$ vertices.

Next we shall discuss the coloring properties of $\mathcal H_{n_1,\ldots,n_s}.$ For a strict coloring $c$ of a mixed hypergraph ${\cal H}=(X, {\cal C},{\cal D})$, we denote by $c(v)$  the color of the vertex $v$, and $c(V)=\{c(v)\mid v\in V\}$.

Note that, for each $i\in [s]$,
$c_i^s=\{X_{i1}^s,X_{i2}^s,\ldots,X_{in_i}^s\}$ is a strict
$n_i$-coloring of ${\cal H}_{n_1,\ldots,n_s}$, where $X_{ij}^s$
consists of the vertices of $X_{n_1,\ldots,n_s}$ whose $i$-th entry is $j.$

\begin{lemma}\label{lem:lem1}  Let $c=\{C_1,C_2,\ldots,C_m\}$ be a strict coloring of ${\cal H}_{n_1,\ldots,n_s}$.
 Then we may reorder the color classes such that the following
 conditions hold:
\begin{itemize}
  \item[\rm(i)] $V_i^s\subseteq C_i$, $U_j^s\subseteq C_j$;
   \item[\rm(ii)] $|c(W_{pk}^s)|=1$, $|c(W_{2h}^s)|=1$\hbox{~and~}  $W_{pk}^s\cap C_i=\emptyset$, $W_{2h}^s\cap C_i=\emptyset$ for any $i\in [r-2]$;
  \item[\rm(iii)] $(n_p+k,\ldots,n_p+k,1,\ldots,1)\notin C_i$ for $i\in [r-2]\setminus
 \{1\}$;
 \item[\rm(iv)]   $(n_1-1,n_2,\ldots,n_{s-1},t,t)\notin C_i$ for $i\in M_s\setminus \{t\}$; and $(n_1,\ldots,n_{s-1},j,j)\notin C_i$ for $i\in M_s\setminus \{j\}$;
  \item[\rm(v)] $(n_s,\ldots,n_s,1,1)\in C_1\cup C_{n_s}$.
 \end{itemize}
\end{lemma}

\proof  (i) Since ${\cal H}_{n_1,\ldots,n_s}[V]$ is complete $r$-uniform,
the $(r-1)^2$ vertices of $V$ fall into exactly $r-1$ color classes, say $C_1,\ldots,C_{r-2}, C_{n_s}$,
 and $|C_i\cap V|=r-1$ for each $i\in M_s$.

We claim that there exists an $i\in M_s$ such that $V_1^s\subseteq
C_i$ or $V_{n_s}^s\subseteq C_i.$
Suppose
 $V_1^s\nsubseteq C_i$ and $V_{n_s}^s\nsubseteq C_i$ for any $i\in M_s$.
Then for each $i\in M_s$, $(V\cap C_i) \cup\{(n_s,\ldots,n_s,1,1)\}$ is a bi-edge, which follows that  $(n_s,\ldots,n_s,1,1)\notin C_i$. Note that for each $i\in M_s$, there
exists an $\alpha_i\in V\cap C_i$ such that
$\{\alpha_1,\ldots,\alpha_{r-2},\alpha_{n_s}\}\neq V_1^s,
V_{n_s}^s$, which imply that  the bi-edge $\{\alpha_1,\ldots,\alpha_{r-2},\alpha_{n_s},(n_s,\ldots,n_s,1,1)\}$ is  polychromatic, a contradiction.
Hence our claim is valid.

Without loss of generality, suppose $V_1^s\subseteq C_1$.  Assume that there exists a $j\in \{2,\ldots,l-1\}$ such that  $|c(V_j^s)|>1$, say  $(j,\ldots,j,1)\in C_2, (j,\ldots,j,2)\in C_3$.  Then  the bi-edge $(V\cap C_i)\cup \{(n_1-1,n_2,\ldots,n_{s-1},j,j)\}$  implies that  $(n_1-1,n_2,\ldots,n_{s-1},j,j)\notin C_i$ for each $i\in  M_s$. It follows that, for any $p\in M_s \setminus \{1,2,3\}$ and $\alpha_p\in C_p\cap V$, the bi-edge $\{(1,\ldots,1),(j,\ldots,j,1),(j,\ldots,j,2),(n_1-1,n_2,\ldots,n_{s-1},j,j),\alpha_4,\ldots,\alpha_{r-2},\alpha_{n_s}\}$   is polychromatic, a contradiction. Hence, $|c(V_i^s)|=1$ for each $i\in \{2,\ldots,l-1\}$. Therefore, we may reorder the color classes such that $V_i^s\subseteq C_i$. Similarly,   we may assume that  $V_i^s\subseteq C_i, i\in \{l,\ldots,r-2\}$. It follows that $V_{n_s}^s\subseteq C_{n_s}$.

Suppose $n_s\geq r$.  Now we focus on the colors of
the vertices of each $U_{j}^s$. For  $i\in M_s$,
the bi-edges $V_i^s\cup \{(j,\ldots,j,t)\},t\in [l]$, imply that
 $U_j^s\cap C_i=\emptyset.$
From the bi-edges $\{(1,\ldots,1),\ldots, (r-2,
\ldots,r-2, 1), (j, \ldots,j, 1), (j, \ldots,j, t)\}, t\in [l]\setminus\{1\}$, we have   $|c(U_j^s)|=1$. Moreover, since ${\cal H}_{n_1,\ldots,n_s}[U_j^s\cup U_q^s]$ is complete $r$-uniform,
$c(U_q^s)\neq c(U_j^s)$ if $q\neq j$. Therefore, we may assume
that $U_j^s\subseteq C_j$. Hence, (i)
holds.

(ii)  From the bi-edges $V_i^s\cup \{(n_p+k,\ldots,n_p+k,n_p,\ldots,n_s,t)\},  t\in [l]\setminus \{1\}$,
we have $W_{pk}^s\cap C_i=\emptyset$ for $i\in [r-2]$. Then the bi-edges
$\{(n_p+k,\ldots,n_p+k,n_p,\ldots,n_s,2),(n_p+k,\ldots,n_p+k,n_p,\ldots,n_s,t),
(1,\ldots,1),\ldots,(r-2,\ldots,r-2,1)\}$, $t\in
[l]\setminus \{1,2\}$,
imply that $|c(W_{pk}^s)|=1$. Similarly, $|c(W_{2h}^s)|=1$ and  $W_{2h}^s\cap C_i=\emptyset$ for any $i\in [r-2]$.
  Hence, (ii) holds.

(iii) The bi-edge
$V_i^s\cup \{(n_p+k,\ldots,n_p+k,1,\ldots,1)\}$
implies that  $(n_p+k,\ldots,n_p+k,1,\ldots,1)\notin C_i$ for $i\in [r-2]\setminus
 \{1\}$.

 (iv) Since $V_i^s\cup\{(n_1-1,n_2,\ldots,n_{s-1},t,t)\}$ is a bi-edge,
 $(n_1-1,n_2,\ldots,n_{s-1},t,t)\notin C_i$ for any  $i\in M_s\setminus \{t\}$; and from the bi-edge $V_i^s\cup\{(n_1,\ldots,n_{s-1},j,j)\}$,  one gets $(n_1,\ldots,n_{s-1},j,j)\notin C_i$ for any $i\in M_s\setminus \{j\}$.

 (v) The bi-edge
$\{(n_s,\ldots,n_s,1,1),(n_s,\ldots,n_s,1),(1,\ldots,1),\ldots,(r-2,\ldots,r-2,1)\}$
implies that  $(n_s,\ldots,n_s,1,1)\in C_1\cup C_{n_s}$. $\qed$

\begin{lemma}\label{lem:lem5}  Let  $c=\{C_1,C_2,\ldots,C_m\}$ be a strict coloring of ${\cal H}_{n_1,\ldots,n_s}$ satisfying the conditions (i)-(v) in Lemma~\ref{lem:lem1}.
 \begin{itemize}
  \item[\rm(i)] If $(n_s,\ldots,n_s,1,1)\in C_1$, then $c=c_s^s$. In particular, $(n_1-1,n_2,\ldots,n_{s-1},t,t)\in C_t$ and $(n_1,\ldots,n_{s-1},j,j)\in C_j$;
  \item[\rm(ii)] If $(n_s,\ldots,n_s,1,1)\in C_{n_s}$, then  $|c(T_1^s)|=1$, $|c(T_2^s)|=1$; and  $T_1^s\cap C_i=\emptyset$,  $T_2^s\cap C_j=\emptyset$ for $i\in M_s$, $j\in [n_s]$.
 \end{itemize}
\end{lemma}

\proof   (i) From the bi-edge
$\{(n_s,\ldots,n_s,1,1),(n_s,\ldots,n_s,2),(2,\ldots,2),\ldots,(r-2,\ldots,r-2,2),(n_1-1,n_2,\ldots,n_{s-1},t,t)\}$, we have $(n_1-1,n_2,\ldots,n_{s-1},t,t)\in C_t$. Similarly, $(n_1,\ldots,n_{s-1},j,j)\in C_j$.

The bi-edges
$\{(n_s,\ldots,n_s,1,1),(2,\ldots,2),\ldots,(r-2,\ldots,r-2,2),(n_s,\ldots,n_s,2),(n_p+k,\ldots,n_p+k,1,\ldots,1)\}$ and $V_{n_s}^s\cup\{(n_p+k,\ldots,n_p+k,1,\ldots,1)\}$ imply that
  $(n_p+k,\ldots,n_p+k,1,\ldots,1)\in C_1$; and  from the bi-edge
$\{(n_s,\ldots,n_s,1,1),(2,\ldots,2),\ldots,(r-2,\ldots,r-2,2),(n_s,\ldots,n_s,t),(n_p+k,\ldots,n_p+k,n_p,\ldots,n_s,t)\}$,
we have  $(n_p+k,\ldots,n_p+k,n_p,\ldots,n_s,t)\in C_{n_s}$. Hence, $c=c_s^s$.

(ii) From the bi-edge $\{(n_s,\ldots,n_s,1),\ldots,(n_s,\ldots,n_s,r-2),(n_s,\ldots,n_s,1,1),(n_1-1,n_2,\ldots,n_{s-1},n_s ,n_s)\}$,  we have $(n_1-1,n_2,\ldots,n_{s-1},n_s,n_s)\notin C_{n_s}$. The bi-edges $\{(n_1-1,n_2,\ldots,n_{s-1},t,t), (n_1-1,n_2,\ldots,n_{s-1},n_s ,n_s),(1,\ldots,1),\ldots,(t-1,\ldots,t-1,1),(t+1,\ldots,t+1,1),\ldots,
(r-2,\ldots,r-2,1),(n_s,\ldots,n_s,1)\}$, $t\in  \{2,\ldots,l-1\}$, imply that $|c(T_1^s)|=1$. Hence,  $T_1^s\cap C_i=\emptyset$ for $i\in M_s$.  Similarly,  $|c(T_2^s)|=1$ and $T_2^s\cap C_i=\emptyset$ for $i\in M_s$. Furthermore, the bi-edge $U_j^s\cup T_2^s$ implies that $T_2^s\cap C_j=\emptyset$ for $j\in [n_s]\setminus M_s$.  Hence, (ii) holds. $\qed$

\begin{lemma}\label{lem:lem7}
 Let  $s\geq 3$ and $c=\{C_1,C_2,\ldots,C_m\}$ be a strict coloring of ${\cal H}_{n_1,\ldots,n_s}$ satisfying the conditions (i)-(iv) in Lemma~\ref{lem:lem1} and $(n_s,\ldots,n_s,1,1)\in C_{n_s}$.
 \begin{itemize}
  \item[\rm(i)] Suppose $c((n_p,\ldots,n_p,1,\ldots,1))\neq c(W_{p0}^s)$ for some $p\in [s-1]\setminus \{1,2\}$. Then   $(n_p+k,\ldots,n_p+k,1,\ldots,1)\in C_1$, and  $W_{pk}^s,T_2^s\subseteq C_a$ for some $a\in [m]\setminus [n_s]$;
  \item[\rm(ii)] Suppose  $c((n_p,\ldots,n_p,1,\ldots,1))=c(W_{p0}^s)$ for some $p\in [s]\setminus \{1,2\}$. Then  $c((n_p+k,\ldots,n_p+k,1,\ldots,1))=c(W_{pk}^s)$;
   \item[\rm(iii)] Suppose $n_1>n_2+1$. If  $c((n_2,1,\ldots,1))\neq c(W_{20}^s)$, then  $(n_2+h,1,\ldots,1)\in C_1$, and  $W_{2h}^s,T_2^s\subseteq C_a$ for some $a\in [m]\setminus [n_s]$; if $c((n_2,1,\ldots,1))=c(W_{20}^s)$, then $c((n_2+h,1,\ldots,1))=c(W_{2h}^s)$.
 \end{itemize}
\end{lemma}

\proof (i) From the bi-edge
$\{(n_p,\ldots,n_p,n_{p+1},\ldots,n_s,2),(1,\ldots,1),\ldots,(r-2,\ldots,$ $r-2,1),(n_p,\ldots,n_p,1,\ldots,1)\}$,
we have  $(n_p,\ldots,n_p,1,\ldots,1)\in C_1$. Then the bi-edge $\{
 (n_p,\ldots,n_p,n_{p+1},\ldots,n_s,2),(n_p,\ldots,n_p,1,\ldots,1),
 (n_1,\ldots,n_{s-1},l,l),(2,\ldots,2),\ldots, (r-2,\ldots,r-2,2)\}$ implies that  $c(W_{p0}^s)=c(T_2^s)$. Suppose
$W_{p0}^s, T_2^s\subseteq C_a$, where $a\in [m]\setminus [n_s]$. From the bi-edge $\{(n_p+k,\ldots,n_p+k,n_p,\ldots,n_s,2),
 (2,\ldots,2),\ldots,(r-2,\ldots,r-2,2),(n_p,\ldots,n_p,n_{p+1},\ldots,n_s,2),(n_p,\ldots,n_p,1,\ldots,1)\}$, one has $W_{pk}^s\subseteq C_a$; and the bi-edges  $W_{pk}^s\cup T_2^s\cup\{(n_p+k,\ldots,n_p+k, 1,\ldots,1)\}$ and $\{(2,\ldots,2),\ldots,(r-2,\ldots,r-2,2),(n_p,\ldots,n_p,n_{p+1},\ldots,n_s,2),
 (n_p,\ldots,n_p,1,\ldots,1),(n_p+k,\ldots,n_p+k,1,\ldots,1)\}$ imply that $(n_p+k,\ldots, n_p+k,1,\ldots,1)\in C_1$.  Hence, (i) holds.

(ii)  If there exists a $k\in [n_{p-1}-n_p-1]$ such that $c((n_p+k,\ldots,n_p+k,1,\ldots,1))\neq c(W_{pk}^s)$, then from the bi-edge
$\{(1,\ldots,1),\ldots,(r-2,\ldots,$ $r-2,1),(n_p+k,\ldots,n_p+k,n_{p},\ldots,n_s,2),(n_p+k,\ldots,n_p+k,1,\ldots,1)\}$, we have
 $(n_p+k,\ldots,n_p+k,1,\ldots,1)\in C_1$. Since $\{
 (n_p+k,\ldots,n_p+k,n_{p},\ldots,n_s,2),(n_p+k,\ldots,n_p+k,1,\ldots,1),(2,\ldots,2),\ldots,(r-2,\ldots,r-2,2),
 (n_1,\ldots,n_{s-1},n_s,n_s)\}$ is a bi-edge,   $W_{pk}^s$ and $T_2^s$
fall into  a common color class, say $C_a$, where $a\in [m]\setminus [n_s]$. The bi-edge $\{(n_p+k,\ldots,n_p+k,n_p,\ldots,n_s,2),(n_p+k,\ldots,n_p+k,1,\ldots,1),
 (n_p,\ldots,n_p,n_{p+1},\ldots,n_s,2),(2,\ldots,2),\ldots,(r-2,\ldots,r-2,2)\}$ implies that $W_{p0}^s\subseteq C_a$; and  the bi-edges  $\{(2,\ldots,2),\ldots,(r-2,\ldots,r-2,2),(n_p+k,\ldots,n_p+k,1,\ldots,1),(n_p,\ldots,n_p,n_{p+1},\ldots,n_s,2),
 (n_p,\ldots,n_p,1,\ldots,1)\}$ and  $W_{p0}^s\cup T_2^s\cup\{(n_p,\ldots,n_p,$ $1,\ldots,1)\}$ imply that $(n_p,\ldots, n_p,1,\ldots,1)\in C_1$, a contradiction. Hence, (ii) holds.

(iii) By the same discussion  in (i) and (ii),  (iii) holds. $\qed$

\begin{lemma}\label{lem:lem2} With the same assumption of  Lemma \ref{lem:lem7}.  Let $b\in [s]\setminus \{1,2\}$ be the minimum number such that $c((n_b,\ldots,n_b,1,\ldots,1))=c(W_{b0}^s)$. Then  for  $p\in \{b,\ldots,s\}$ and $k\in [n_{p-1}-n_p-1]\cup \{0\}$,  we may reorder the color classes such that the following conditions hold:
\begin{itemize}
  \item[\rm(i)] $\{(n_p+k,\ldots,n_p+k,1,\ldots,1)\}\cup W_{pk}^s\subseteq C_{n_p+k}$;
  \item[\rm(ii)] if $b>3$, then $c=c_{b-1}^s$.
\end{itemize}
\end{lemma}

\proof The bi-edge $\{(n_s,\ldots,n_s,1,1),(n_s,\ldots,n_s,2),\ldots,(n_s,\ldots,n_s,r-1),(n_p+k,\ldots,n_p+k,n_p,\ldots,n_s,2)\}$ implies that $W_{pk}^s\cap C_{n_s}=\emptyset$.

 (i)  We claim that   $c((n_p,\ldots,n_p,1,\ldots,1))=c(W_{p0}^s)$. Suppose  there exists a $p\in  \{b+1,\ldots,s\}$ such that $c((n_p,\ldots,n_p,1,\ldots,1))\neq c(W_{p0}^s)$. Then by Lemma~\ref{lem:lem7} we have $(n_p,\ldots,n_p,1,\ldots,1)\in C_1$, and $W_{p0}^s, T_2^s\subseteq C_a$ for some $a\in [m]\setminus [n_s]$. The bi-edge $\{(n_p,\ldots,n_p,1,\ldots,1),
 (n_p,\ldots,n_p,n_{p+1},\ldots,n_s,2),(2,\ldots,2),\ldots,(r-2,\ldots,r-2,2),(n_b,\ldots,n_b,n_{b+1},\ldots,n_s,2)\}$ implies that $W_{b0}^s\subseteq C_a$; and from the bi-edges $W_{b0}^s\cup T_2^s\cup\{(n_b,\ldots,n_b,$ $1,\ldots,1)\}$ and $\{
 (n_b,\ldots,n_b,n_{b+1},\ldots,n_s,2),(2,\ldots,2),\ldots,(r-2,\ldots,r-2,2),(n_p,\ldots,n_p,1,\ldots,1),(n_b,\ldots,n_b,1,\ldots,1)
 \}$, we have $(n_b,\ldots, n_b,1,\ldots,1)\in C_1$, a contradiction. Hence, our claim is valid.

  By Lemma~\ref{lem:lem7}, we have  $c((n_p+k,\ldots,n_p+k,1,\ldots,1))=c(W_{pk}^s)$. Furthermore, since  ${\cal H}_{n_1,\ldots,n_s}[\{(n_{p_1}+k_1,\ldots,n_{p_1}+k_1,1,\ldots,1),(n_{p_2}+k_2,\ldots,n_{p_2}+k_2,1,\ldots,1)\}\cup W_{p_1k_1}^s\cup W_{p_2k_2}^s]$ is complete $r$-uniform,  $c(W_{p_1k_1}^s)\neq c(W_{p_2k_2}^s)$ for $p_1,p_2\in \{b,\ldots,s\}$ and $k_1\in [n_{p_1-1}-n_{p_1}-1]\cup \{0\}, k_2\in [n_{p_2-1}-n_{p_2}-1]\cup \{0\}$. For any $j\in [n_s]\setminus M_s$, the fact that ${\cal H}_{n_1,\ldots,n_s}[U_j^s\cup W_{pk}^s\cup \{(n_p+k,\ldots,n_p+k,1,\ldots,1)\}]$ is complete $r$-uniform implies that $W_{pk}^s\cap C_j=\emptyset$. Hence, we may reorder the color classes such  that $\{(n_p+k,\ldots,n_p+k,1,\ldots,1)\}\cup W_{pk}^s\subseteq C_{n_p+k}$, as desired.

  (ii) Suppose $b>3$. That is to say, $c((n_q,\ldots,n_q,1,\ldots,1))\neq c(W_{q0}^s)$ for each $q\in \{3,\ldots,b-1\}$. Then by Lemma~\ref{lem:lem7} we have $(n_q+j,\ldots,n_q+j,1,\ldots,1)\in C_1$, and  $W_{qj}^s, T_2^s\subseteq C_a$ for some $a\in [m]\setminus [n_s]$. Since $\{
 (2,\ldots,2),\ldots,(r-2,\ldots,r-2,2),(n_2+h,n_2,\ldots,n_s,2),(n_3,n_3,n_3,n_4,\ldots,n_s,2),(n_3,n_3,1,\ldots,1)\}$ is a bi-edge, $W_{2h}^s\subseteq C_a$; and from the bi-edges $\{(n_2+h,1,\ldots,1),
 (n_3,n_3,n_3,n_4,\ldots,n_s,2),(2,\ldots,2),\ldots,(r-2,\ldots,r-2,2),(n_3,n_3,1,\ldots,1)\}$ and $\{(n_2+h,1,\ldots,1)\}\cup W_{2h}^s\cup T_2^s$, we have $(n_2+h,1,\ldots,1)\in C_1$. Similarly, $(n_1-1,1,\ldots,1)\in C_1$ and $T_1^s\subseteq C_a$.  Moreover, the bi-edge $\{(n_p+k,\ldots,n_p+k,1,\ldots,1)\}\cup W_{pk}^s\cup T_2^s$ implies that $T_2^s\cap C_{n_p+k}=\emptyset$. Hence,  $W_{qj}^s, W_{2h}^s, T_1^s, T_2^s \subseteq C_{n_{b-1}}$, which implies that $c=c_{b-1}^s$.  $\qed$

\section{The proof of Theorem~\ref{thm}}

In this section we shall
 construct  $r$-uniform
bi-hypergraphs which are one-realizations of $S$ and meet  the bound
in Lemma~\ref{lem:ineq} in each case. We divide our constructions into
three subsections.

\subsection{ $n_1>n_2+1$ or $n_1=n_2+1>r+2$}

\begin{thm}\label{thm:thm20}  ${\cal H}_{n_1,n_2}$ is a one-realization of $\{n_1,n_2\}$.
\end{thm}

\proof It suffices
to prove that
 $c_1^2, c_2^2$ are all the strict
colorings of ${\cal H}_{n_1,n_2}$.
Suppose $c=\{C_1,C_2,\ldots,C_m\}$ is a strict coloring of
${\cal H}_{n_1,n_2}$ satisfying  the conditions (i)-(v) in   Lemma~\ref{lem:lem1}. In particular, $(n_2,1,1)\in
C_1\cup C_{n_2}.$

\textbf{Case 1}\,\ $(n_2,1,1)\in C_1$.

By Lemma~\ref{lem:lem5} we have $c=c_2^2$.

  \textbf{Case 2} \,\ $(n_2,1,1)\in C_{n_2}$.

From the bi-edges $\{(n_2,1,1),(n_2,n_2,2),\ldots,(n_2,n_2,r-1),(n_2+h,n_2,2)\}$ and $\{(n_2+h,1,1),(n_2,n_2,1),\ldots,(n_2,n_2,r-1)\}$, one gets $W_{2h}^2\cap C_{n_2}=\emptyset$ and $(n_2+h,1,1)\notin C_{n_2}$. The bi-edge $\{(2,2,2),\ldots,(r-2,r-2,2),(n_2,n_2,2),(n_2+h,1,1),(n_2+h,n_2,2)\}$ implies that $c((n_2+h,1,1))=c(W_{2h}^2)$. Since ${\cal H}_{n_1,n_2}[\{(n_2+h,1,1)\}\cup W_{2h}^2\cup U_j^2]$ is complete $r$-uniform,  $W_{2h}^2\cap C_j=\emptyset$ for any  $j\in [n_2]\setminus M_2$. Note that ${\cal H}_{n_1,n_2}[\{(n_2+h_1,1,1),(n_2+h_2,1,1)\}\cup W_{2h_1}^2\cup W_{2h_2}^2]$ is complete $r$-uniform.  Hence, we may reorder the color classes such  that $\{(n_2+h,1,1)\}\cup W_{2h}^2\subseteq C_{n_2+h}$. By Lemma~\ref{lem:lem5}, one gets that  $|c(T_1^2)|=|c(T_2^2)|=1$;  $T_1^2\cap C_i=\emptyset$ and $T_2^2\cap C_j=\emptyset$ for any $i\in M_2, j\in [n_2]$.  From the bi-edges $\{(n_1-1,1,1),(n_1-1,2,2),(n_2,n_2,2),(2,2,2),\ldots,
(r-2,r-2,2)\}$ and $\{(n_2,n_2,1),\ldots,(n_2,n_2,r-1),(n_1-1,1,1)\}$, we have  $c((n_1-1,1,1))=c(T_1^2)$. Then the fact that  ${\cal H}_{n_1,n_2}[\{(n_1-1,1,1)\}\cup T_1^2\cup U_j^2]$ is complete $r$-uniform implies that $T_1^2\cap C_j=\emptyset$ for any $j\in [n_2]\setminus M_2$.  And since ${\cal H}_{n_1,n_2}[\{(n_1-1,1,1),(n_2+h,1,1)\}\cup T_1^2\cup W_{2h}^2]$ is complete $r$-uniform, $T_1^2\cap C_{n_2+h}=\emptyset$. Therefore, $T_1^2\subseteq C_{n_1-1}$.
The bi-edge $\{(n_2+h,1,1)\}\cup W_{2h}^2\cup T_2^2$ implies that $T_2^2\cap C_{n_2+h}=\emptyset$. Moreover, from the bi-edge $\{(n_1-1,1,1)\}\cup T_1^2\cup T_2^2$, we have $T_2^2\cap C_{n_1-1}=\emptyset$.
Hence, $c=c_1^2$.$\qed$

 \begin{thm} \label{thm:thm1} Suppose $s\geq 3$.  Then  ${\cal H}_{n_1,\ldots,n_s} $ is a one-realization of
$\{n_1,n_2,\ldots,n_s\}$.
\end{thm}

\proof It suffices
to prove that
 $c_1^s,\ldots, c_s^s$ are all the strict
colorings of ${\cal H}_{n_1,\ldots,n_s}$.
Suppose $c=\{C_1,C_2,\ldots,C_m\}$ is a strict coloring of
${\cal H}_{n_1,\ldots,n_s}$ satisfying  the conditions (i)-(v) in   Lemma~\ref{lem:lem1}. In particular, $(n_s,\ldots,n_s,1,1)\in
C_1\cup C_{n_s}.$

\textbf{Case 1}\,\ $(n_s,\ldots,n_s,1,1)\in C_1$.

By Lemma~\ref{lem:lem5}  we have $c=c_s^s$.

  \textbf{Case 2} \,\ $(n_s,\ldots,n_s,1,1)\in C_{n_s}$.

Then $c$ satisfies the condition (ii) in Lemma~\ref{lem:lem5}.
In this case, we shall prove that $c\in \{c_1^s,\ldots, c_{s-1}^s\}$.  Let $b\in [s]\setminus \{1,2\}$ be the minimum number such that $c((n_b,\ldots,n_b,1,\ldots,1))=c(W_{b0}^s)$.  By Lemma~\ref{lem:lem2},  we have $\{(n_p+k,\ldots,n_p+k,1,\ldots,1)\}\cup W_{pk}^s\subseteq C_{n_p+k}$ for
 $p\in \{b,\ldots,s\}$ and $k\in
[n_{p-1}-n_p-1]\cup \{0\}$.

  \textbf{Case 2.1} \,\  $b>3$. By Lemma~\ref{lem:lem2}  we have  that $c=c_{b-1}^s$.

  \textbf{Case 2.2} \,\  $b=3$.   Suppose $n_1>n_2+1$. We  focus on the colors of the vertices of  $\{(n_2,1,\ldots,1)\}\cup W_{20}^s$.  If $c((n_2,1,\ldots,1))\neq c(W_{20}^s)$, then  by Lemma~\ref{lem:lem7} we have  $(n_2+h,1,\ldots,1)\in C_1$, and $W_{2h}^s, T_2^s\subseteq C_a$ for some $a\in [m]\setminus [n_s]$. The bi-edge $\{(n_2,1,\ldots,1),
 (n_2,n_2,n_3,n_4,\ldots,n_s,2),(2,\ldots,2),\ldots,(r-2,\ldots,r-2,2),(n_1-1,n_2,\ldots,n_{s-1},2,2)\}$ implies that $T_1^s\subseteq C_a$. From the bi-edges $\{(n_1-1,1,\ldots,1)\}\cup T_1^s\cup T_2^s$ and
$\{(n_2,1,\ldots,1),
 (n_2,n_2,n_3,n_4,\ldots,n_s,2),(n_1-1,1,\ldots,1),(2,\ldots,2),\ldots,(r-2,\ldots,r-2,2)\}$, we have $(n_1-1,1,\ldots,1)\in C_1$; and the bi-edge $\{(n_p+k,\ldots,n_p+k,1,\ldots,1)\}\cup W_{pk}^s\cup T_2^s$ implies that $T_2^s\cap C_{n_p+k}=\emptyset$. Hence,  $W_{2h}^s, T_1^s, T_2^s \subseteq C_{n_2}$, which imply that $c=c_2^s$.
 If $c((n_2,1,\ldots,1))=c(W_{20}^s)$,
then by  Lemma~\ref{lem:lem7} we have that $c((n_2+h,1,\ldots,1))=c(W_{2h}^s)$. Since  ${\cal H}_{n_1,\ldots,n_s}[\{(n_2+h_1,1,\ldots,1), (n_2+h_2,1,\ldots,1)\}\cup W_{2h_1}^s\cup W_{2h_2}^s]$ is complete $r$-uniform,  $c(W_{2h_1}^s)\neq c(W_{2h_2}^s)$ if $h_1\neq h_2$. The fact that ${\cal H}_{n_1,\ldots,n_s}[\{(n_p+k,\ldots,n_p+k,1,\ldots,1),(n_2+h,1,\ldots,1)\}\cup W_{pk}^s\cup W_{2h}^s]$ is complete $r$-uniform implies that $W_{2h}^s\cap C_{n_p+k}=\emptyset$. Moreover, for each $j\in [n_s]\setminus M_s$, since ${\cal H}_{n_1,\ldots,n_s}[U_j^s\cup W_{2h}^s\cup \{(n_2+h,1,\ldots,1)\}]$ is complete $r$-uniform, $W_{2h}^s\cap C_j=\emptyset$.  Hence,  $\{(n_2+h,1,\ldots,1)\}\cup W_{2h}^s\subseteq C_{n_2+h}$. From the bi-edges $\{(n_1-1,1,\ldots,1),(n_1-1,n_2,\ldots,n_{s-1},2,2),(1,\ldots,1),\ldots,(r-2,\ldots,r-2,1)\}$ and $\{(n_1-1,1,\ldots,1),(n_1-1,n_2,\ldots,n_{s-1},2,2),(n_s,\ldots,n_s,2),(2,\ldots,2),\ldots,(r-2,\ldots,r-2,2)\}$, Similarly, $\{(n_1-1,1,\ldots,1)\}\cup T_1^s\subseteq C_{n_1-1}$ and $T_2^s\subseteq C_{n_1}$. It follows that  $c=c_1^s$.

 Suppose $n_1=n_2+1$.  If $c(T_1^s)=c(T_2^s)$, then from the bi-edge $\{(n_2,1,\ldots,1)\}\cup T_1^s\cup T_2^s$, we have  $c((n_2,1,\ldots,1))\neq c(T_1^s)$. The bi-edge $\{(n_2,1,\ldots,1),(1,\ldots,1),\ldots,(r-2,\ldots,r-2,1),(n_2,n_2,n_3,\ldots,n_{s-1},2,2)\}$ implies that $(n_2,1,\ldots,1)\in C_1$. The bi-edge $\{(n_p+k,\ldots,n_p+k,1,\ldots,1)\}\cup W_{pk}^s\cup T_2^s$ implies that $T_2^s\cap C_{n_p+k}=\emptyset$. Hence, $c=c_2^s$. If $c(T_1^s)\neq c(T_2^s)$,
  then from the bi-edges $\{(n_2,n_2,n_3,\ldots,n_{s-1},2,2),(1,\ldots,1),\ldots,(r-2,\ldots,r-2,1),(n_2,1,\ldots,1)\}$ and  $\{(n_1,\ldots,n_{s-1},2,2),(n_2,1,\ldots,1),(2,\ldots,2),\ldots,(r-2,\ldots,r-2,2),(n_2,n_2,n_3,\ldots,n_{s-1},2,2)\}$, we have $c((n_2,1,\ldots,1))=c(T_1^s)$. Since   ${\cal H}_{n_1,\ldots,n_s}[\{(n_p+k,\ldots,n_p+k,1,\ldots,1)\}\cup W_{pk}^s\cup T_1^s\cup \{(n_2,1,\ldots,1)\}]$ is complete $r$-uniform,  $T_1^s\cap C_{n_p+k}=\emptyset$; and the bi-edge $T_2^s\cup W_{pk}^s\cup \{(n_p+k,\ldots,n_p+k,1,\ldots,1)\}$ implies that $T_2^s\cap C_{n_p+k}=\emptyset$.  Moreover, for each $j\in [n_s]\setminus M_s$, the fact that ${\cal H}_{n_1,\ldots,n_s}[U_j^s\cup T_1^s\cup \{(n_1-1,1,\ldots,1)\}]$ is complete $r$-uniform implies that $T_1^s\cap C_j=\emptyset$. Therefore,  $c=c_1^s$.  $\qed$

  Note that ${\cal H}_{n_1,\ldots,n_s}$ is a desired one-realization of $S$ when $n_1>n_2+1$, or $s\geq 3$, $n_1=n_2+1>r+2$, $r$ is even.

   Moreover, from the proof of Theorem~\ref{thm:thm1} we observe  that if $s\geq 3$, $n_1=n_2+1$ and $r$ is odd, then for any $j\in [n_s]\setminus M_s$, $p\in [s]\setminus \{1,2\}$ and $k\in [n_{p-1}-n_p-1]\cup \{0\}$, $U_j^s\cup T_1^s$ and $T_1^s\cup W_{pk}^s\cup \{(n_p+k,\ldots,n_p+k,1,\ldots,1)\}$ are bi-edges,  so   ${\cal H}_{n_1,\ldots,n_s}[X']$ is also a  one-realization of $S$, where $X'=X_{n_1,\ldots,n_s}\setminus \{(n_2,1,\ldots,1)\}$.

 \begin{thm} \label{thm:thm14} Let $X'=X_{n_1,\ldots,n_s}\setminus \{(n_2,1,\ldots,1)\}$.  Then  ${\cal H}_{n_1,\ldots,n_s}[X']$ is a  desired one-realization of $S$ when $s\geq 3$, $n_1=n_2+1>r+2$ and $r$ is odd.
\end{thm}

 \medskip
{\bf Construction \uppercase\expandafter{\romannumeral2}.}  Let $l=\lceil\frac r2\rceil$. Write
\begin{eqnarray*}
V_i&=&\{(i,i,t)\mid t\in [r-1]\}, i\in [r-1]; ~ U_j=\{(j,j,k)\mid k\in [l]\}, j=r,\ldots,n_2,\\
W&=&\{(n_1,n_2,1),\ldots,(n_1,n_2,l)\},\\
Y&=&(\bigcup\limits_{i=1}^{r-1}V_i)\cup (\bigcup\limits_{j=r}^{n_2}U_j)\cup W\cup \{(r-1,1,1)\},\\
{\cal E}&=&\{\{\alpha_1,\ldots,\alpha_r\}\mid\alpha_q\in Y,   2\leq |\{\alpha_1^j,\ldots,\alpha_r^j\}|\leq r-1,q\in [r], j\in [3]\}\\
&\cup& \{\{(r-1,1,1),(1,1,1),\ldots,(r-1,r-1,1)\}\}\\
&\cup&\{\{(r-1,1,1),(1,1,1),\ldots,(1,1,r-1)\}\}.\end{eqnarray*}

\begin{lemma}\label{lem:lem6}  Suppose $n_1=n_2+1> r+2$. Let $e=\{E_1,E_2,\ldots,E_m\}$ be a strict coloring of ${\cal G}$. Then we may reorder the color classes such that the following conditions hold:

\begin{itemize}
  \item[\rm(i)] $V_i\subseteq E_i$, $U_j\subseteq E_j$;
  \item[\rm(ii)] $(r-1,1,1)\in E_1\cup E_{r-1}$;
  \item[\rm(iii)] $|e(W)|=1$ and  $W\cap E_i=\emptyset$ for any $i\in [n_2-1]$.
 \end{itemize}
\end{lemma}

\proof (i) Let $T=\cup_{i=1}^{r-1}V_i$. Referring to the proof of Lemma~\ref{lem:lem1}, we have that
the $(r-1)^2$ vertices of $T$ fall into exactly $r-1$ color classes,
say $E_1,E_2,\ldots, E_{r-1}$, and $|E_i\cap T|=r-1$ for each $i\in [r-1]$. Moreover,  there exists an $i\in [r-1]$  such that $V_1\subseteq
E_i$ or $V_{r-1}\subseteq E_i.$
 Without loss generality, suppose $V_1\subseteq E_1$. Then, $\alpha\notin E_1$ for any $\alpha\in T\setminus V_1$. Note that $n_2>r+1$. For any $i\in [r-1]$, since $(T\cap E_i)\cup \{(r,r,1)\}$ is a bi-edge, $(r,r,1)\notin E_i$.   Suppose $(2,2,1)\in E_2$. If $V_2\cap E_2\neq V_2$, say $(2,2,2)\in E_3$, then pick
   $\alpha_i\in T\cap E_i, i\in \{4,\ldots,r-1\}$, the bi-edge  $\{(1,1,1),(2,2,1),(2,2,2),(r,r,1),\alpha_4,\ldots,\alpha_{r-1}\}$
  is polychromatic, a contradiction.  Hence,  $|e(V_2)|=1$, say $V_2\subseteq E_2$. Similarly, we may reorder the color classes such that  $V_i\subseteq E_i$ for $i\in [r-1]$.

  Now we focus on the colors of
the vertices of each $U_{j}$. For  $i\in [r-1]$,
the bi-edges $V_i\cup \{(j,j,k)\},k\in [l]$, imply that
 $U_j\cap E_i=\emptyset.$
From the bi-edges $\{(1,1,1),\ldots, (r-2,
r-2, 1), (j, j, 1), (j,j,k)\}, k\in [l]\setminus\{1\}$, we have
$|e(U_j)|=1$. Moreover, since ${\cal G}[U_q\cup U_j]$ is complete $r$-uniform,
 $e(U_q)\neq e(U_j)$ if $q\neq j$. Therefore, we may assume
that $U_j\subseteq E_j$.

(ii) For $i\in [r-2]\setminus \{1\}$, the bi-edge $\{(r-1,1,1),(i,i,1),\ldots,(i,i,r-1)\}$ implies that $(r-1,1,1)\notin E_i$. Then from the bi-edge $\{(r-1,1,1),(1,1,1),\ldots,(r-1,r-1,1)\}$, we have $(r-1,1,1)\in E_1\cup E_{r-1}$.

(iii) For $t\in [l]$ and $i\in [r-1]$, since $\{(n_1,n_2,t),(i,i,1),\ldots,(i,i,r-1)\}$ is a bi-edge,
$W\cap E_i=\emptyset$. From the bi-edges $\{(n_1,n_2,1),(n_1,n_2,t),(1,1,1),\ldots,(r-2,r-2,1)\}, t\in [l]\setminus \{1\}$, we have $|e(W)|=1$. For any $j\in \{r,\ldots,n_2-1\}$, the fact that ${\cal G}[W\cup U_j]$ is complete $r$-uniform implies that $W\cap E_j=\emptyset$. Hence, (iii) holds. $\qed$

\begin{thm}\label{thm:thm3}  Let $n_1=n_2+1> r+2$. Then ${\cal G}$ is a one-realization of
$\{n_1,n_2\}$.
\end{thm}

\proof Note that
\begin{eqnarray*}&&e_1=\{V_1,\ldots,V_{r-2},V_{r-1}\cup \{(r-1,1,1)\}, U_r,\ldots, U_{n_2-1}, U_{n_2}, W\},\\
 &&e_2=\{V_1,\ldots,V_{r-2},V_{r-1}\cup \{(r-1,1,1)\}, U_r,\ldots, U_{n_2-1}, U_{n_2}\cup W\}\end{eqnarray*} are two strict colorings of ${\cal G}$.
  Let  $e=\{E_1,E_2,\ldots,E_m\}$ be a strict coloring of ${\cal G}$ satisfying the conditions (i)-(iii) in Lemma~\ref{lem:lem6}. Then the bi-edge  $\{(1,1,1),\ldots,(1,1,r-1),(r-1,1,1)\}$ implies that $(r-1,1,1)\in E_{r-1}$. Hence,
$e=e_2$ or $e=e_1$ if $e(U_{n_2})=e(W)$ or not. $\qed$

Note that,  when  $n_1=n_2+1> r+2$ and $r$ is even,  $|Y|=(n_1-3)r-(n_1-r-3)\lfloor\frac{r}{2}\rfloor+2$, which implies that ${\cal G}$ is a desired one-realization of
$\{n_1,n_2\}$.

\begin{thm}\label{thm:thm15}  Let  $Y'=Y\setminus \{(n_2,n_2,1),(n_1,n_2,1)\}$. Then ${\cal G}'={\cal G}[Y']$ is a  desired one-realization of $\{n_1,n_2\}$ when  $n_1=n_2+1>r+2$ and $r$ is odd.
\end{thm}

\proof Note that 
\begin{eqnarray*}&&e_1'=\{V_1,\ldots,V_{r-2},V_{r-1}\cup \{(r-1,1,1)\}, U_r,\ldots, U_{n_2-1}, U_{n_2}', W'\},\\
 &&e_2'=\{V_1,\ldots,V_{r-2},V_{r-1}\cup \{(r-1,1,1)\}, U_r,\ldots, U_{n_2-1}, U_{n_2}'\cup W'\}\end{eqnarray*} are two strict colorings of ${\cal G}$, where $U_{n_2}'=U_{n_2}\setminus \{(n_2,n_2,1)\}$ and $W'=W\setminus \{(n_1,n_2,1)\}$.
  Let  $e'=\{E_1',E_2',\ldots,E_m'\}$ be a strict coloring of ${\cal G}'$. Then by Lemma~\ref{lem:lem6}, we may reorder the color classes such that
  \begin{itemize}
  \item[\rm(i)]  $V_i\subseteq E_i'$,  $U_j\subseteq E_j'$;
  \item[\rm(ii)] $|e'(W')|=1$ and   $W'\cap E_i'=\emptyset$ for any $i\in [r-1]$;
  \item[\rm(iii)] $(r-1,1,1)\in E_1'\cup E_{r-1}'$.
 \end{itemize}
Similarly, we have $|e'(U_{n_2}')|=1$ and  $U_{n_2}'\cap E_i'=\emptyset$ for any $i\in [r-1]$.
 Moreover,  since $U_{n_2}'\cup U_j$ and $U_j\cup W'$ are bi-edges, $U_{n_2}'\cap E_j'=\emptyset$ and $W'\cap E_j'=\emptyset$. Hence, $e'=e_2'$ or  $e'=e_1'$ if $e'(U_{n_2}')=e'(W')$ or not.  
 
Since $|Y'|=(n_1-3)r-(n_1-r-3)\lfloor\frac{r}{2}\rfloor+2$, the desired result follows.$\qed$

\subsection{$n_1=n_2+1=r+2$}

Let $n_1=n_2+1=r+2$.  There are the
following four possible subcases: $S=\{r+2,r+1\}$, $S=\{r+2,r+1,r\}$, $S=\{r+2,r+1,r-1\}$ or $S=\{r+2,r+1,r,r-1\}$.

\medskip
{\bf Construction \uppercase\expandafter{\romannumeral3}.}  Let
\begin{eqnarray*}
X_1&=&\{(i,i,t)\mid i\in [r],t\in [r-1]\}\cup \{(r-1,1,1),(r+1,r+1,2),(r+2,r+1,2)\},\\
{\cal B}_1&=&\{\{\alpha_1,\ldots,\alpha_r\}\mid\alpha_q\in X_1,   2\leq |\{\alpha_1^j,\ldots,\alpha_r^j\}|\leq r-1, q\in [r], j\in [3]\}\\
&\cup& \{\{(r-1,1,1),(1,1,1),\ldots,(r-1,r-1,1)\}\}\\
&\cup&\{\{(r-1,1,1),(1,1,1),\ldots,(1,1,r-1)\}\}.\end{eqnarray*}
 Then ${\cal H}_1=(X_1, {\cal
B}_1)$ is an $r$-uniform bi-hypergraph with $r(r-1)+3$ vertices.

\begin{thm} \label{thm:thm5} ${\cal H}_1$ is a desired one-realization of
$\{r+2,r+1\}$.
\end{thm}

\proof For $i\in [r]$, let $X_{1i}=\{(i,i,t)|~t\in [r-1]\}$ and $Z_{1,r-1}=X_{1,r-1}\cup \{(r-1,1,1)\}$.
Then
\begin{eqnarray*}&&c_{11}=\{X_{11},\ldots,X_{1,r-2},Z_{1,r-1}, X_{1r},\{(r+1,r+1,2)\},\{(r+2,r+1,2)\}\},\\
 &&c_{12}=\{X_{11},\ldots,X_{1,r-2},Z_{1,r-1}, X_{1r},\{(r+1,r+1,2),(r+2,r+1,2)\}\}\end{eqnarray*} are two
strict colorings of ${\cal H}_1$.
  For any strict coloring $c=\{C_1,C_2,\ldots,C_m\}$  of ${\cal H}_1$, by  Lemma~\ref{lem:lem6} we may reorder the color classes and get that
\begin{itemize}
  \item[\rm(i)] $X_{1i}\subseteq C_i$ for $i\in [r]$;
  \item[\rm(ii)] $(r+1,r+1,2),(r+2,r+1,2)\notin C_i$ for
$i\in [r]$;
  \item[\rm(iii)] $(r-1,1,1)\in C_1\cup C_{r-1}$.
 \end{itemize}
Since
$\{(r-1,1,1),(1,1,1),\ldots,(1,1,r-1)\}$ is a bi-edge,  $(r-1,1,1)\in C_{r-1}$. Therefore, $c=c_{12}$ or  $c=c_{11}$ if
$c((r+1,r+1,2))=c((r+2,r+1,2))$ or not. $\qed$

     \medskip
{\bf Construction \uppercase\expandafter{\romannumeral4}.}  Let
\begin{eqnarray*}
X_2&=&\{(i,i,i,t)\mid i\in [r],t\in [r-1]\}\cup \{(r,r,1,1),(n_2,n_2,r,2),(n_1,n_2,r,2)\},\\
{\cal B}_2&=&\{\{\alpha_1,\ldots,\alpha_r\}\mid \alpha_q\in X_2,   2\leq |\{\alpha_1^j,\ldots,\alpha_r^j\}|\leq r-1, q\in [r], j\in [4]\}\\
&\cup& \{\{(r,r,1,1),(1,1,1,1),\ldots,(r-2,r-2,r-2,1),(r,r,r,1)\}\}.\end{eqnarray*}
 Then ${\cal H}_2=(X_2, {\cal
B}_2)$ is an $r$-uniform bi-hypergraph with $r(r-1)+3$ vertices.

\begin{thm}\label{thm:thm6} Suppose $n_1=n_2+1=r+2$ and $n_3=r$. Then ${\cal H}_2$ is a desired one-realization of
$\{n_1,n_2,n_3\}$.
\end{thm}

\proof For  $i\in [r]$, let $X_{2i}=\{(i,i,i,t)|~t\in [r-1]\}$.
Note that
\begin{eqnarray*}&&c_{21}=\{X_{21},\ldots,X_{2,r-1}, X_{2r}\cup \{(r,r,1,1)\}, \{(n_2,n_2,r,2)\},\{(n_1,n_2,r,2)\}\},\\
&&c_{22}=\{X_{21},\ldots,X_{2,r-1}, X_{2r}\cup \{(r,r,1,1)\}, \{(n_2,n_2,r,2),(n_1,n_2,r,2)\}\},\\
 &&c_{23}=\{X_{21}\cup \{(r,r,1,1)\},X_{22},\ldots,X_{2,r-1}, X_{2r}\cup \{(n_2,n_2,r,2),(n_1,n_2,r,2)\}\}\end{eqnarray*} are three
strict colorings of ${\cal H}_2$.
  For any strict coloring $c=\{C_1,C_2,\ldots,C_m\}$ of  ${\cal H}_2$, by Lemma~\ref{lem:lem6}
 we may reorder the color classes such  that $X_{2i}\subseteq C_i$ for $i\in [r]$; $(n_2,n_2,r,2),(n_1,n_2,r,2)\notin C_j$ for $j\in [r-1]$; and $(r,r,1,1)\in C_1\cup C_r$.

 \textbf{Case 1}\,\  $(r,r,1,1)\in C_1$.

 The two bi-edges
$\{(3,3,3,1),\ldots,(r,r,r,1),(r,r,1,1),(n_2,n_2,r,2)\}$ and
 $\{(3,3,3,1), $ $ \ldots,(r,r,r,1),(r,r,1,1),(n_1,n_2,r,2)\}$
  imply that $(n_2,n_2,r,2), (n_1,n_2,r,2)\in C_r$.
  Hence, $c=c_{23}$.

\textbf{Case 2}\,\ $(r,r,1,1)\in C_r$.

The two bi-edges
$\{(r,r,1,1),(r,r,r,2),\ldots,(r,r,r,r-1),(n_2,n_2,r,2)\}$ and
 $\{(r,r,1,1),$ $(r,r,r,2),\ldots,(r,r,r,r-1),(n_1,n_2,r,2)\}$
 imply that $(n_2,n_2,r,2),(n_1,n_2,r,2)\notin C_r$.
 Hence, $c=c_{22}$ or  $c=c_{21}$ if $c((n_2,n_2,r,2))=c((n_1,n_2,r,2))$ or not. $\qed$

    \medskip
{\bf Construction \uppercase\expandafter{\romannumeral5}.}  Let
\begin{eqnarray*}
X_3&=&\{(i,i,i,t)\mid i,t\in [r-1]\}\cup \{(r,r,t,t)\mid t\in [r-1]\}\\
&\cup& \{(r-1,r-1,1,1),(n_2,n_2,r-1,2),(n_1,n_2,r-1,2)\},\\
{\cal B}_3&=&\{\{\alpha_1,\ldots,\alpha_r\}\mid \alpha_q\in X_3,   2\leq |\{\alpha_1^j,\ldots,\alpha_r^j\}|\leq r-1, q\in [r], j\in [4]\}\\
&\cup& \{\{(r-1,r-1,1,1),(1,1,1,1),\ldots,(r-1,r-1,r-1,1)\}\}.\end{eqnarray*}
 Then ${\cal H}_3=(X_3, {\cal
B}_3)$ is an $r$-uniform bi-hypergraph with $r(r-1)+3$ vertices.

\begin{thm}\label{thm:thm7}  Suppose $n_1=n_2+1=r+2, n_3=r-1$. Then ${\cal H}_3$ is a desired one-realization of
$\{n_1,n_2,n_3\}$.
\end{thm}

\proof
Note that for any $i\in [3]$, $c_{3i}=\{X_{3i}^1,X_{3i}^2,\ldots,X_{3i}^{n_i}\}$ is a strict $n_i$-coloring of ${\cal H}_3$, where $X_{3i}^j$  consists of the vertices of $X_3$ whose $i$-th entry is $j$.

   For any strict coloring $c=\{C_1,C_2,\ldots,C_m\}$ of ${\cal H}_3$, by Lemma~\ref{lem:lem1} we may assume that   $(i,i,i,t)\in C_i$;  $(r-1,r-1,1,1)\in C_1\cup C_{r-1}$; and for any $i\in [r-1]\setminus \{j\}$, $(r,r,j,j)\notin C_i$. Moreover,  the bi-edges $\{(i,i,i,1),\ldots,(i,i,i,r-1),(n_2,n_2,r-1,2)\}$ and $\{(i,i,i,1),\ldots,(i,i,i,r-1),(n_1,n_2,r-1,2)\}$ imply that $(n_2,n_2,r-1,2),(n_1,n_2,r-1,2)\notin C_i$ for any $i\in [r-2]$.

  \textbf{Case 1}\,\ $(r-1,r-1,1,1)\in C_1$.

  For $j\in [r-1]\setminus \{1\}$,  the bi-edges $\{(r-1,r-1,1,1),(2,2,2,j),\ldots,(r-1,r-1,r-1,j),(r,r,j,j)\}$ imply  that $(r,r,j,j)\in C_j$. From the bi-edge $\{(r,r,1,1),(r-1,r-1,1,1),(2,2,2,2),\ldots,(r-1,r-1,r-1,2)\}$, we have $(r,r,1,1)\in C_1$. Since  $\{(r,r,1,1),\ldots,(r,r,r-1,r-1),(n_2,n_2,r-1,2)\}$ and $\{(r,r,1,1),\ldots,(r,r,r-1,r-1),(n_1,n_2,r-1,2)\}$
  are bi-edges,  $(n_2,n_2,r-1,2), (n_1,n_2,r-1,2)\in C_{r-1}$. It follows that $c=c_{33}$.

   \textbf{Case 2}\,\ $(r-1,r-1,1,1)\in C_{r-1}$.

   The bi-edges $\{(n_2,n_2,r-1,2),(r-1,r-1,1,1),(r-1,r-1,r-1,2),\ldots,(r-1,r-1,r-1,r-1)\}$ and $\{(n_1,n_2,r-1,2),(r-1,r-1,1,1),(r-1,r-1,r-1,2),\ldots,(r-1,r-1,r-1,r-1)\}$ imply that $(n_2,n_2,r-1,2), (n_1,n_2,r-1,2)\notin C_{r-1}$. Note that $(r,r,1,1)\notin C_i$ for any $i\in [r-1]\setminus \{1\}$. Suppose $(r,r,1,1)\in C_1$. Then for $j\in [r-1]\setminus \{1\}$, the bi-edge $\{(r,r,1,1),(r,r,j,j),(2,2,2,2),\ldots,(j-1,j-1,j-1,2),(j+1,j+1,j+1,2),\ldots,
   (r-2,r-2,r-2,2),(r-1,r-1,1,1),(n_2,n_2,r-1,2)\}$  implies that $c((r,r,j,j))=c((n_2,n_2,r-1,2))$. Similarly, $c((r,r,j,j))=c((n_1,n_2,r-1,2))$. Therefore, the bi-edge $\{(r,r,2,2),\ldots,(r,r,r-1,r-1),(n_2,n_2,r-1,2),(n_1,n_2,r-1,2)\}$ is monochromatic, a contradiction.
  Hence, $(r,r,1,1)\in C_r$.
   From the bi-edge $\{(1,1,1,2),(3,3,3,2),\ldots,(r-2,r-2,r-2,2),(r-1,r-1,1,1),(r,r,1,1),(r,r,2,2)\}$, we have $(r,r,2,2)\in C_r$. Similarly, $(r,r,j,j)\in C_r$ for any $j\in [r-1]$. Since $\{(n_2,n_2,r-1,2),(r,r,1,1),\ldots,(r,r,r-1,r-1)\}$ and $\{(n_1,n_2,r-1,2),(r,r,1,1),\ldots,(r,r,r-1,r-1)\}$ are bi-edges, $(n_2,n_2,r-1,2), (n_1,n_2,r-1,2)\notin  C_r$. Therefore, $c=c_{32}$ or $c=c_{31}$ if $c((n_2,n_2,r-1,2))=c((n_1,n_2,r-1,2))$ or not. $\qed$

\medskip
{\bf Construction \uppercase\expandafter{\romannumeral6}.}  Let
\begin{eqnarray*}
X_4&=&\{(i,i,i,i,t)\mid i,t\in [r-1]\}\cup \{(r,r,r,t,t)\mid t\in [r-1]\setminus \{1\}\}\\
&\cup&  \{(r-1,r-1,r-1,1,1),(r,r,1,1,1),(n_2,n_2,r,r-1,2),(n_1,n_2,r,r-1,2)\},\\
{\cal B}_4&=&\{\{\alpha_1,\ldots,\alpha_r\}\mid\alpha_q\in X_4,   2\leq |\{\alpha_1^j,\ldots,\alpha_r^j\}|\leq r-1,q\in [r], j\in [5]\}\\
&\cup& \{\{(r-1,r-1,r-1,1,1),(1,1,1,1,1),\ldots,(r-1,r-1,r-1,r-1,1)\}\}.\end{eqnarray*}
 Then ${\cal H}_4=(X_4, {\cal
B}_4)$ is an $r$-uniform bi-hypergraph with $r(r-1)+3$ vertices.

\begin{thm}\label{thm:thm8} Suppose $n_1=n_4+3=r+2$. Then ${\cal H}_4$ is a desired one-realization of
$\{n_1,n_2,n_3,n_4\}$.
\end{thm}

\proof Note that, for any $i\in [4]$, $c_{4i}=\{X_{4i}^1,X_{4i}^2,\ldots,X_{4i}^{n_i}\}$ is a strict $n_i$-coloring of ${\cal H}_4$, where $X_{4i}^j$  consists of the vertices of $X_4$ whose $i$-th entry is $j$.

  For any strict coloring $c=\{C_1,C_2,\ldots,C_m\}$ of ${\cal H}_4$, by Theorem~\ref{thm:thm7} we may assume that   $(i,i,i,i,t)\in C_i$; $(r,r,r,t,t)\notin C_j$ for any  $j\in [r-1]\setminus \{t\}$; $(r,r,1,1,1)\notin C_i$ for $i\in [r-1]\setminus \{1\}$;  $(n_2,n_2,r,r-1,2),(n_1,n_2,r,r-1,2)\notin C_j$ for  $j\in [r-2]$; and $(r-1,r-1,r-1,1,1)\in C_1\cup C_{r-1}$.

  \textbf{Case 1}\,\ $(r-1,r-1,r-1,1,1)\in C_1$.

  Referring to the proof of Theorem~\ref{thm:thm7}, we have  $c=c_{44}$.

   \textbf{Case 2}\,\ $(r-1,r-1,r-1,1,1)\in C_{r-1}$.

  The bi-edges $\{(n_2,n_2,r,r-1,2),(r-1,r-1,r-1,1,1),(r-1,r-1,r-1,r-1,2),\ldots,(r-1,r-1,r-1,r-1,r-1)\}$ and $\{(n_1,n_2,r,r-1,2),(r-1,r-1,r-1,1,1),(r-1,r-1,r-1,r-1,2),\ldots,(r-1,r-1,r-1,r-1,r-1)\}$ imply that $(n_2,n_2,r,r-1,2), (n_1,n_2,r,r-1,2)\notin C_{r-1}$.

 \textbf{Case 2.1}\,\ $(r,r,1,1,1)\in C_1$. The bi-edge
  $\{(r-1,r-1,r-1,1,1),(n_2,n_2,r,r-1,2),(r,r,r,t,t),(r,r,1,1,1),(2,2,2,2,2),
  \ldots,(t-1,t-1,t-1,t-1,2),(t+1,t+1,t+1,t+1,2),\ldots,(r-2,r-2,r-2,r-2,2)\}$  implies that $c((r,r,r,t,t))=c((n_2,n_2,r,r-1,2))$. Similarly, $c((r,r,r,t,t))=c((n_1,n_2,r,r-1,2))$. Hence, $c=c_{43}$.

   \textbf{Case 2.2}\,\  $(r,r,1,1,1)\notin C_1$. Then we may assume that $(r,r,1,1,1)\in C_r$.
   From the bi-edge $\{(1,1,1,1,2),(3,3,3,3,2),\ldots,(r-2,r-2,r-2,r-2,2),(r-1,r-1,r-1,1,1),(r,r,1,1,1),(r,r,r,2,2)\}$, we have $(r,r,r,2,2)\in C_r$. Similarly, $(r,r,r,t,t)\in C_r$ for any $t\in [r-1]\setminus \{1\}$. Since $\{(n_2,n_2,r,r-1,2),(r,r,1,1,1),(r,r,r,2,2),\ldots,(r,r,r,r-1,r-1)\}$ and $\{(n_1,n_2,r,r-1,2),(r,r,1,1,1),(r,r,r,2,2),\ldots,(r,r,r,r-1,r-1)\}$ are bi-edges, $(n_2,n_2,r,r-1,2), (n_1,n_2,r,r-1,2)\notin  C_r$. Hence, $c=c_{42}$ or $c=c_{41}$ if $c((n_2,n_2,r-1,2))=c((n_1,n_2,r-1,2))$ or not.
 $\qed$

\subsection{ $n_1=n_2+1\leq r+1$}

Let $n_1=n_2+1\leq r+1$.  There are the
following three  possible subcases: $S=\{r+1,r\}$, $S=\{r+1,r,r-1\}$ or $S=\{r,r-1\}$.

\medskip
{\bf Construction \uppercase\expandafter{\romannumeral7}.} Let
\begin{eqnarray*}
X_5&=&\{(i,i,t)\mid i,t\in  [r-1]\}\cup \{(r-1,1,1),(r,r,2),(r+1,r,2)\},\\
{\cal B}_5&=&\{\{\alpha_1,\ldots,\alpha_r\}\mid\alpha_q\in X_5,   2\leq |\{\alpha_1^j,\ldots,\alpha_r^j\}|\leq r-1, q\in [r], j\in [3]\}\\
&\cup& \{\{(r-1,1,1),(1,1,1),\ldots,(r-1,r-1,1)\}\}\\
&\cup&\{\{(r-1,1,1),(1,1,1),\ldots,(1,1,r-1)\}\}.\end{eqnarray*}
 Then ${\cal H}_5=(X_5, {\cal
B}_5)$ is an $r$-uniform bi-hypergraph with $(r-1)^2+3$ vertices.

\begin{thm}\label{thm:thm9}  ${\cal H}_5$ is a desired one-realization of
$\{r+1,r\}$.
\end{thm}

\proof For any $i\in [r-1]$, let $X_{5i}=\{(i,i,t)|~t\in [r-1]\}$.
Note that
\begin{eqnarray*}&&c_{51}=\{X_{51},\ldots,X_{5,r-2},X_{5,r-1}\cup \{(r-1,1,1)\}, \{(r,r,2)\},\{(r+1,r,2)\}\},\\
 &&c_{52}=\{X_{51},\ldots,X_{5,r-2},X_{5,r-1}\cup \{(r-1,1,1)\}, \{(r,r,2),(r+1,r,2)\}\}\end{eqnarray*} are two
strict colorings of ${\cal H}_5$.
  For any strict coloring $c=\{C_1,C_2,\ldots,C_m\}$  of ${\cal H}_5$, referring to the proof of  Lemma~\ref{lem:lem6},  we may assume that  $X_{5i}\subseteq C_i$; $(r,r,2),(r+1,r,2)\notin C_i$ for
$i\in [r-1]$; $(r-1,1,1)\in C_1\cup C_{r-1}$. Since
$\{(r-1,1,1),(1,1,1),\ldots,(1,1,r-1)\}$ is a bi-edge,  $(r-1,1,1)\in C_{r-1}$. Therefore, $c=c_{52}$ or $c=c_{51}$ if
$c((r,r,2))=c((r+1,r,2))$ or not.
$\qed$

\medskip
{\bf Construction \uppercase\expandafter{\romannumeral8}.}  Let
\begin{eqnarray*}
X_6&=&\{(i,i,i,t)\mid i,t\in [r-1]\}\cup \{(r,r,j,j)\mid j=3,\ldots,r-2\}\\
&\cup &\{(r-1,r-1,1,1), (r+1,r,r-1,r-1)\},\\
{\cal B}_6&=&\{\{\alpha_1,\ldots,\alpha_r\}\mid\alpha_q\in X_6,   2\leq |\{\alpha_1^j,\ldots,\alpha_r^j\}|\leq r-1, q\in [r], j\in [4]\}\\
&\cup& \{\{(r-1,r-1,1,1),(1,1,1,1),\ldots,(r-1,r-1,r-1,1)\}\}.\end{eqnarray*}
 Then ${\cal H}_6=(X_6, {\cal
B}_6)$ is an $r$-uniform bi-hypergraph with $r(r-1)-1$ vertices.

\begin{thm}\label{thm:thm10} Suppose $r\geq 6$. Then ${\cal H}_6$ is a desired one-realization of
$\{r+1,r,r-1\}$.
\end{thm}

\proof For $i\in [3]$, $c_{6i}=\{X_{6i}^1,\ldots,X_{6i}^{n_i}\}$ is a strict $n_i$-coloring of ${\cal H}_6$, where $X_{6i}^j$ consists of the vertices of $X_6$ whose $i$-th entry is $j$.

Let  $X_{6i}=\{(i,i,i,t)\mid t\in [r-1]\}, i=1,2,\ldots,r-1$, and $V=\cup_{i=1}^{r-1}X_{6i}$. For a strict coloring
 $c=\{C_1,C_2,\ldots,C_m\}$  of ${\cal H}_6$,
referring to the proof of Lemma~\ref{lem:lem1},  we have that
the $(r-1)^2$ vertices of $V$ fall into exactly $r-1$ color classes,
say $C_1,C_2,\ldots, C_{r-1}$,  each of which  contains exactly
$r-1$ vertices of $V$. Furthermore, there exists an $i\in [r-1]$ such that $X_{61}\subseteq
C_i$ or $X_{6,r-1}\subseteq C_i.$
Without loss of generality, suppose $X_{61}\subseteq C_1$.  Assume that there exists a $k\in \{3,\ldots,r-2\}$ such that  $|c(X_{6k})|>1$, say $k=3$ and $(3,3,3,1)\in C_2, (3,3,3,2)\in C_3$.  Then the bi-edge $(V\cap C_i)\cup \{(r,r,3,3)\}$ implies that  $(r,r,3,3)\notin C_i$ for any $i\in [r-1]$. Hence, for $j\in \{4,\ldots,r-1\}$ and $\alpha_j\in C_j\cap V$, the bi-edge $\{(1,1,1,1),(3,3,3,1),(3,3,3,2),(r,r,3,3),\alpha_4,\ldots,\alpha_{r-1}\}$ is polychromatic,  a contradiction. It follows that $|c(X_{6i})|=1$ for any $i\in \{3,\ldots,r-2\}$.  Therefore, we may assume that $X_{6i}\subseteq C_i$ for each $i\in \{3,\ldots,r-2\}$. If  $|c(X_{6,r-1})|>1$, say $(r-1,r-1,r-1,1)\in C_2,(r-1,r-1,r-1,2)\in C_{r-1}$, then for each $i\in [r-1]$, from the bi-edge $(C_i\cap V)\cup \{(r+1,r,r-1,r-1)\}$, one gets $(r+1,r,r-1,r-1)\notin C_i$. It follows that the bi-edge $\{(1,1,1,1),(3,3,3,1),\ldots,(r-2,r-2,r-2,1),(r-1,r-1,r-1,1),(r-1,r-1,r-1,2),(r+1,r,r-1,r-1)\}$ is polychromatic, a contradiction. So $|c(X_{6,r-1})|=1$, which follows that  $|c(X_{62})|=1$.
Hence, we may assume that $X_{6i}\in C_i, i\in [r-1]$.  For $i\in [r-1]\setminus \{j\}$, the bi-edge $\{(r,r,j,j),(i,i,i,1),\ldots,(i,i,i,r-1)\}$ implies that $(r,r,j,j)\notin C_i$. From the bi-edge $\{(r-1,r-1,1,1),(1,1,1,1),\ldots,(r-1,r-1,r-1,1)\}$, one has
 $(r-1,r-1,1,1)\in C_1\cup C_{r-1}$.

 \textbf{Case 1 }\,\  $(r-1,r-1,1,1)\in C_1$.

 From the bi-edge $\{(r-1,r-1,1,1),(2,2,2,1),\ldots,(r-1,r-1,r-1,1),(r,r,j,j)\}$, we have $(r,r,j,j)\in C_j$. Similarly, $(r+1,r,r-1,r-1)\in C_{r-1}$. It follows that $c=c_{63}$.

  \textbf{Case 2 }\,\  $(r-1,r-1,1,1)\in C_{r-1}$.

 From the bi-edge $\{(r+1,r,r-1,r-1),(r-1,r-1,1,1),(r-1,r-1,r-1,2),\ldots,(r-1,r-1,r-1,r-1)\}$, we have $(r+1,r,r-1,r-1)\notin C_{r-1}$. Note that $(r,r,3,3)\notin C_i$ for $i\in [r-1]\setminus \{3\}$. If $(r,r,3,3)\in C_3$, then the bi-edge $\{(1,1,1,1),(2,2,2,1),(5,5,5,1),\ldots,(r-1,r-1,r-1,1),(r,r,3,3),(r,r,4,4),(r+1,r,r-1,r-1)\}$ implies that $c((r,r,4,4))=c((r+1,r,r-1,r-1))$. It follows that $(r,r,4,4)\notin C_4$, and the bi-edge $\{(4,4,4,1),\ldots,(r-1,r-1,r-1,1),(r,r,3,3),(r,r,4,4),(1,1,1,1),(2,2,2,1)\}$ is polychromatic, a contradiction. Suppose $(r,r,3,3)\in C_r$. For $j\in \{4,\ldots,r-2\}$, since $\{(1,1,1,1),\ldots,(j-1,j-1,j-1,1),(j+1,j+1,j+1,1),\ldots,(r-1,r-1,r-1,1),(r,r,3,3),(r,r,j,j)\}$ is a bi-edge,  $(r,r,j,j)\in C_r$. Hence,  $c=c_{62}$ if $(r+1,r,r-1,r-1)\in  C_r$;  $c=c_{61}$ if  $(r+1,r,r-1,r-1)\notin C_r$. $\qed$

 \medskip
{\bf Construction \uppercase\expandafter{\romannumeral9}.}  Let
\begin{eqnarray*}
X_7&=&\{(1,1,1,t),(3,3,3,t)\mid t\in [4]\}\cup \{(2,2,2,j)\mid j\in [3]\}\\
&\cup& \{(4,4,4,1),(4,4,4,3),(4,4,4,4),(4,4,2,2),(2,2,4,4)\}\\
&\cup &\{(4,4,1,1),(5,5,2,2),(6,5,4,4)\},\\
{\cal B}_7&=&\{\{\alpha_1,\alpha_2,\alpha_3,\alpha_4,\alpha_5\}\mid\alpha_q\in X_7,   2\leq |\{\alpha_1^j,\alpha_2^j,\alpha_3^j,\alpha_4^j,\alpha_5^j\}|\leq 4, q\in [5], j\in [4]\}\\
&\cup& \{\{(4,4,1,1),(1,1,1,1),(2,2,2,1),(3,3,3,1),(4,4,4,1)\}\}.\end{eqnarray*}
 Then ${\cal H}_7=(X_7, {\cal
B}_7)$ is an $r$-uniform bi-hypergraph with $r(r-1)-1$ vertices.

\begin{thm}\label{thm:thm11}  ${\cal H}_7$ is a desired one-realization of
$\{6,5,4\}$.
\end{thm}

\proof For $i\in [3]$, $c_{7i}=\{X_{7i}^1,X_{7i}^2,\ldots,X_{7i}^{n_i}\}$ is a strict $n_i$-coloring of ${\cal H}_7$, where $X_{7i}^j$ consists of the vertices of $X_7$ whose $i$-entry is $j$.
Let \begin{eqnarray*}
W_1&=&\{(1,1,1,t)\mid t\in [4]\},~~~~~~W_2=\{(2,2,2,j)\mid j\in [3]\},\\
W_3&=&\{(3,3,3,t)\mid t\in [4]\},~~~~~~W_4=\{(4,4,4,1),(4,4,4,3),(4,4,4,4),(4,4,2,2)\},\\
W_4'&=&\{(4,4,4,k)\mid k=1,3,4\},~W=W_1\cup W_2\cup W_3\cup W_4\cup \{(2,2,4,4)\}.
\end{eqnarray*}
 For any strict coloring $c=\{C_1,C_2,\ldots,C_m\}$ of ${\cal H}_7$, since ${\cal H}_7[W]$ is complete $r$-uniform,  the vertices of $W$ fall into 4 color classes, say $C_1,C_2,C_3,C_4$, and $|W\cap C_i|=4$, $i\in [4]$.  We claim that there exists an $i\in [4]$ such that $W_1\subseteq C_i$ or $W_4\subseteq C_i$.  Suppose $W_1\nsubseteq C_i, W_4\nsubseteq C_i$ for any $i\in [4]$.  Then  the bi-edge $(W\cap C_i)\cup \{(4,4,1,1)\}$  implies that $(4,4,1,1)\notin C_i$ for $i\in [4]$. Assume that $(4,4,1,1)\in C_5$. Note that for each $i\in [4]$, there exists an $\alpha_i\in W\cap C_i$ such that $\{\alpha_1,\alpha_2,\alpha_3,\alpha_4,(4,4,1,1)\}$ is a bi-edge. It follows that the bi-edge $\{\alpha_1,\alpha_2,\alpha_3,\alpha_4,(4,4,1,1)\}$  is polychromatic, a contradiction. Hence, our claim is valid.

\textbf{Case 1 }\,\   $|c(W_1)|=1$.

Let $W_1\subseteq C_1$. If $|c(W_4')|>1$, then  the bi-edge $(C_i\cap W)\cup \{(6,5,4,4)\}$ implies that $(6,5,4,4)\notin C_i$ for each $i\in [4]$. Note that for any $i\in [4]$, there exists a $\beta_i\in W\cap C_i$ such that the bi-edge $\{\beta_1,\beta_2,\beta_3,\beta_4,(6,5,4,4)\}$  is polychromatic, a contradiction. Hence, $|c(W_4')|=1$, say $W_4'\subseteq C_4$. Similarly,  $|c(W_2)|=1$, say $W_2\subseteq C_2$.  If $(2,2,4,4)\notin C_2\cup C_4$, then $(2,2,4,4)\in C_3$. Hence,  from the bi-edge $(C_i\cap W)\cup \{(6,5,4,4)\}$, we have  $(6,5,4,4)\notin C_i$ for any $i\in [4]$.  Note that for any $i\in [4]$, there exists a $\gamma_i\in W\cap C_i$ such that the bi-edge $\{\gamma_1,\gamma_2,\gamma_3,\gamma_4,(6,5,4,4)\}$  is polychromatic, a contradiction. Therefore, $(2,2,4,4)\in C_2\cup C_4$. Similarly, $(4,4,2,2)\in C_2\cup C_4$. So $W_3\subseteq C_3$. Moreover, from the bi-edges $\{(4,4,1,1),(1,1,1,1),(2,2,2,1),(3,3,3,1),(4,4,4,1)\}$, $\{(4,4,1,1)\}\cup W_3$ and  $\{(4,4,1,1)\}\cup (W\cap C_2)$, one gets $(4,4,1,1)\in C_1\cup C_4$.

\textbf{Case 1.1 }\,\  $(2,2,4,4)\in C_2, (4,4,2,2)\in C_4$. Then $(5,5,2,2),(6,5,4,4)\notin C_i, i\in [4]$. Since $\{(4,4,1,1),(2,2,4,4),(3,3,3,1),(4,4,4,1),(6,5,4,4)\}$ is a bi-edge,  $(4,4,1,1)\in C_4$. Hence, $c=c_{72}$ or $c=c_{71}$ if $c((5,5,2,2))=c((6,5,4,4))$ or not.

\textbf{Case 1.2 }\,\ $(2,2,4,4)\in C_4, (4,4,2,2)\in C_2$. Then we have $(4,4,1,1)\notin C_4$ from the bi-edge $\{(4,4,1,1),(4,4,4,1),(4,4,4,3),(4,4,4,4),(2,2,4,4)\}$, which follows that $(4,4,1,1)\in C_1$. Note that $(5,5,2,2)\notin C_1\cup C_3\cup C_4$ and $(6,5,4,4)\notin C_1\cup C_2\cup C_3$. The bi-edge $\{(4,4,1,1),(2,2,2,2),(3,3,3,1),(4,4,4,1),(5,5,2,2)\}$ implies that $(5,5,2,2)\in C_2$; and since $\{(6,5,4,4),(1,1,1,1),(2,2,2,1),(3,3,3,1),(2,2,4,4)\}$ is a bi-edge, $(6,5,4,4)\in C_4$. Hence $c=c_{73}$.

\textbf{Case 2 }\,\   $|c(W_4)|=1$.

Let $W_4\subseteq C_4$. Then  the bi-edge $(W\cap C_i)\cup \{(5,5,2,2)\}$ implies that $(5,5,2,2)\notin C_i$ for each $i\in [4]$. If $|c(W_1)|>1$, say $(1,1,1,1)\in C_1, (1,1,1,2)\in C_2$, then for any $\alpha_3\in (W\cap C_3)$, the bi-edge $\{(1,1,1,1),(1,1,1,2),\alpha_3, (4,4,4,1),(5,5,2,2)\}$ is polychromatic, a contradiction. Hence $|c(W_1)|=1$. Similarly,  $|c(W_3)|=1$, which follows that  $|c(W_2\cup \{(2,2,4,4)\})|=1$. Without loss of generality, suppose $W_1\subseteq C_1, W_3\subseteq C_3, W_2\cup \{(2,2,4,4)\}\subseteq C_2$. By Case 1.1 we  have $(4,4,1,1)\in C_4$. Hence, $c=c_{72}$ or $c=c_{71}$ if $c((5,5,2,2))=c((6,5,4,4))$ or not. $\qed$

 \medskip
{\bf Construction \uppercase\expandafter{\romannumeral10}.}  Let
\begin{eqnarray*}
X_8&=&\{(1,1,1,t)\mid t\in [3]\}\cup \{(2,2,2,1),(2,2,2,2)\}\cup \{(3,3,3,1),(3,3,3,3)\}\\
&\cup &\{(3,3,1,1),(2,2,3,3),(3,3,2,2),(4,4,2,2),(5,4,3,3)\},\\
{\cal B}_8&=&\{\{\alpha_1,\alpha_2,\alpha_3,\alpha_4\}\mid\alpha_q\in X_8,   2\leq |\{\alpha_1^j,\alpha_2^j,\alpha_3^j,\alpha_4^j\}|\leq 3, q\in [4], j\in [4]\}\\
&\cup& \{\{(3,3,1,1),(1,1,1,1),(2,2,2,1),(3,3,3,1)\}\}.\end{eqnarray*}
 Then ${\cal H}_8=(X_8, {\cal
B}_8)$ is an $r$-uniform bi-hypergraph  with $(r-1)^2+3$ vertices. Furthermore, similar to Theorem~\ref{thm:thm11}, we have the following result.

\begin{thm}\label{thm:thm12}  ${\cal H}_8$ is a desired one-realization of
$\{5,4,3\}$.
\end{thm}

\medskip
{\bf Construction \uppercase\expandafter{\romannumeral11}.}  Let
\begin{eqnarray*}
X_9&=&\{(i,i,t)\mid i,t\in [r-1]\}\cup \{(r,j,j)\mid j=3,\ldots,r-1\}\cup \{(r-1,1,1)\},\\
{\cal B}_9&=&\{\{\alpha_1,\ldots,\alpha_r\}\mid\alpha_q\in X_9,   2\leq |\{\alpha_1^j,\ldots,\alpha_r^j\}|\leq r-1, q\in [r], j\in [3]\}\\
&\cup& \{\{(r-1,1,1),(1,1,1),\ldots,(r-1,r-1,1)\}\}.\end{eqnarray*}
 Then ${\cal H}_9=(X_9, {\cal
B}_9)$ is an $r$-uniform bi-hypergraph with $r(r-1)-1$ vertices.

Similar to Theorem~\ref{thm:thm10}, we have the following result.

\begin{thm}\label{thm:thm13} Suppose $r\geq 4$.  Then ${\cal H}_9$ is a desired one-realization of
$\{r,r-1\}$.
\end{thm}

Combining  Lemmas~\ref{lem:ineq1},\ref{lem:ineq} and   Theorems~\ref{thm:thm20}-\ref{thm:thm14},  Theorems~\ref{thm:thm3}-\ref{thm:thm13},
the proof of Theorem 1.1 is completed.

\section*{Acknowledgment}

 This research is partially supported by NSF of
Shandong Province (No. ZR2009AM013), NSF of China
(No. 11226288),  NSFC(11271047), and  the
Fundamental Research Funds for the Central University of China and Troy
University Research Grant.

\end{CJK*}

\end{document}